\documentclass[reqno]{amsart}
\usepackage{hyperref}
\usepackage{amsmath}
\usepackage[backend=bibtex,style=alphabetic,firstinits,useprefix=true,url=false]{biblatex}
\usepackage{amssymb}
\usepackage{amsthm}
\usepackage{mathrsfs}
\usepackage{tikz}
\usetikzlibrary{arrows,matrix}

\numberwithin{equation}{section}
\newtheorem{theorem}{Theorem}[section]
\newtheorem*{claim}{\textit{Claim}}
\newtheorem{proposition}[theorem]{Proposition}
\newtheorem{corollary}[theorem]{Corollary}
\newtheorem{lemma}[theorem]{Lemma}
\newtheorem{conjecture}[theorem]{Conjecture}

\newcommand{\define}[1]{\emph{#1}}
\newcommand{\nbar}{\vert\!\vert}
\newcommand{\bnbar}{\Bigg\vert\!\Bigg\vert}
\renewcommand{\epsilon}{\varepsilon}
\newcommand{\intd}{\,\mathrm{d}}
\newcommand{\condex}[2]{\mathbb{E}({#1}\vert{#2})}
\newcommand{\haar}{\mathrm{m}}

\DeclareMathOperator{\ap}{\mathcal{A}}
\DeclareMathOperator{\eig}{\mathcal{E}}
\DeclareMathOperator{\wm}{\mathcal{W}}
\DeclareMathOperator{\conv}{\star}
\DeclareMathOperator{\fp}{FP}
\DeclareMathOperator{\ip}{IP}
\DeclareMathOperator{\lp}{L}
\DeclareMathOperator{\ret}{R}
\DeclareMathOperator{\total}{\ast}
\DeclareMathOperator{\sep}{Sep}
\DeclareMathOperator{\symdiff}{\triangle}
\DeclareMathOperator{\dens}{d}
\DeclareMathOperator{\upperdens}{\overline{d}}
\DeclareMathOperator{\lowerdens}{\underline{d}}

\begin{document}

\begin{abstract}
We describe characteristic factors for certain averages arising from commuting actions of locally compact, second-countable, amenable groups.
Under some ergodicity assumptions we use these factors to prove a form of multiple recurrence for three such actions.
\end{abstract}

\title{Characteristic Factors for Commuting Actions of Amenable Groups}
\author{Donald Robertson}
\address{Department of Mathematics\\
The Ohio State University\\
231 West 18th Avenue\\
Columbus\\
OH 43210-1174\\
USA}
\email{robertson@math.ohio-state.edu}
\date{\today}

\maketitle

\section{Introduction}
\label{sec:introduction}

Furstenberg and Katznelson's multiple recurrence theorem \cite{MR531279} states that if $T_1,\dots,T_k$ are commuting, measure-preserving transformations of a probability space $(X,\mathscr{B},\mu)$ then
\begin{equation*}
\liminf_{N \to \infty} \frac{1}{N} \sum_{n=1}^N \mu(B \cap T_1^{-n} B \cap \cdots \cap T_k^{-n} B) > 0
\end{equation*}
for any $B$ in $\mathscr{B}$ with $\mu(B) > 0$.
It is natural to ask whether such a result holds for commuting actions of groups, by which we mean actions $T_1,\dots,T_k$ of a group $G$ on a probability space $(X,\mathscr{B},\mu)$ by measure-preserving transformations that satisfy $T_i^g T_j^h = T_j^h T_i^g$ for all $g,h \in G$ and all $1 \le i < j \le k$.
Unfortunately the results in \cite{MR1180251} suggest that in certain cases
\begin{equation*}
\mu(B \cap (T_1^g)^{-1} B \cap \cdots \cap (T_k^g)^{-1} B) = 0
\end{equation*}
for all $g \ne 1$ in $G$.
However, if one instead considers multiple recurrence of the form
\begin{equation*}
\mu(B \cap (T_k^g \cdots T_1^g)^{-1} B \cap (T_k^g \cdots T_2^g)^{-1} B \cap \cdots \cap (T_k^g)^{-1} B) > 0
\end{equation*}
then the situation is more promising.
Bergelson, McCutcheon and Zhang \cite{MR1481813} proved that when $G$ is countable and amenable and $\mu(B) > 0$ the set
\begin{equation}
\label{eqn:amenableRothReturnTimes}
\{ g \in G \,:\, \mu(B \cap (T_2^g T_1^g)^{-1} B \cap (T_2^g)^{-1} B) > 0 \}
\end{equation}
is syndetic, meaning that finitely many of its left-shifts cover $G$.
In fact, Bergelson and McCutcheon \cite{MR2354320} have shown for any countable group $G$ that \eqref{eqn:amenableRothReturnTimes} belongs to any minimal idempotent ultrafilter in $\beta G$.
Also, it follows from the work of Bergelson and Rosenblatt (Theorem~2.4 in \cite{MR961735}) that if $G$ is amenable and if $T_j \cdots T_i$ is weakly-mixing for all $1 \le i \le j \le k$ then
\begin{equation*}
\{ g \in G \,:\, \mu(B \cap (T_k^g \cdots T_1^g)^{-1} B \cap \cdots \cap (T_k^g)^{-1} B) \ge \mu(B)^{k+1} \}
\end{equation*}
has full density with respect to any F{\o}lner sequence in $G$.
More generally, Bergelson has made the following conjecture.

\begin{conjecture}[Section 5 of \cite{MR1411215}]
Let $G$ be a countable amenable group with a left F\o{}lner sequence $\Phi$.
Let $T_1,\dots,T_k$ be commuting, measure-preserving actions of $G$ on a probability space $(X,\mathscr{B},\mu)$.
Then
\begin{equation}
\label{eqn:amenableMultipleRecurrence}
\liminf_{N \to \infty} \frac{1}{|\Phi_N|} \sum_{g \in \Phi_N} \mu(B \cap (T_k^g \cdots T_1^g)^{-1} B \cap \cdots \cap (T_k^g)^{-1}B) > 0
\end{equation}
for any $B \in \mathscr{B}$ with $\mu(B) > 0$.
\end{conjecture}

In this paper we describe characteristic factors for the average
\begin{equation}
\label{eqn:amenableCommutingAverage}
\frac{1}{|\Phi_N|} \sum_{g \in \Phi_N} \prod_{i=1}^k T_k^g \cdots T_i^g f_i
\end{equation}
that will allow us to verify a version of Bergelson's conjecture when $k = 3$ under the assumption that the actions $T_1,T_2$ and $T_2T_1$ are ergodic.
(In fact we will do so for locally-compact, second-countable, amenable groups, but only discuss the discrete case in the introduction.)
By \define{characteristic factors} we mean $T_k \cdots T_i$ invariant sub-$\sigma$-algebras $\mathscr{C}_{k,i}$ of $\mathscr{B}$ such that
\begin{equation*}
\lim_{N \to \infty} \frac{1}{|\Phi_N|} \sum_{g \in \Phi_N} \bigg( \prod_{i=1}^k T_k^g \cdots T_i^g f_i - \prod_{i=1}^k T_k^g \cdots T_i^g \condex{f_i}{\mathscr{C}_{k,i}} \bigg) = 0
\end{equation*}
in $\lp^2(X,\mathscr{B},\mu)$ for any $f_i$ in $\lp^\infty(X,\mathscr{B},\mu)$.
By identifying characteristic factors we reduce the study of the limiting behavoir of \eqref{eqn:amenableCommutingAverage} to the situation where $f_i$ is $\mathscr{C}_{k,i}$ measurable.

Although the terminology is more recent, this technique was first used by Furstenberg in his ergodic proof \cite{MR0498471} of Szemeredi's theorem.
Therein he exhibited, for any ergodic, measure-preserving transformation $T$ of a probability space $(X,\mathscr{B},\mu)$, an increasing sequence $\mathscr{Z}_k$ of $T$ invariant sub-$\sigma$-algebras, with $\mathscr{Z}_k$ an isometric extension of $\mathscr{Z}_{k-1}$, such that
\begin{equation*}
\lim_{N-M \to \infty} \frac{1}{N-M} \sum_{n=M}^{N-1} \int \prod_{i=1}^k T^{in} f_i \cdot f_{k+1} - \prod_{i=1}^k T^{in} \condex{f_i}{\mathscr{Z}_k} \cdot f_{k+1} \intd\mu = 0
\end{equation*}
in $\lp^2(X,\mathscr{B},\mu)$ for any $f_i$ in $\lp^\infty(X,\mathscr{B},\mu)$.
Furstenberg then used the properties of isometric extensions to show by induction on $i$ that
\begin{equation}
\label{eqn:SzPositivity}
\liminf_{N-M \to \infty} \frac{1}{N-M} \sum_{n=M}^{N-1} \mu(B \cap T^{-n}B \cap \cdots \cap T^{-kn} B) > 0
\end{equation}
for any $B$ in $\mathscr{Z}_i$ having positive measure.

More recently, Host and Kra~\cite{MR2150389} and Ziegler~\cite{MR2257397} have shown that one can replace $\mathscr{Z}_k$ with a smaller sub-$\sigma$-algebra that corresponds to an inverse limit of $k$-step nilrotations.
This has lead to sharper (e.g. \cite{MR2138068}, \cite{MR2435427}) combinatorial results.
Also, Frantzikinakis and Kra~\cite{MR2142946} have shown, under natural ergodicity assumptions, that inverse limits of commuting rotations on a nilmanifold are characteristic for commuting $\mathbb{Z}$ actions.

Our techniques are similar to those used in \cite{MR0498471}.
However, since we deal with commuting actions, our characteristic factors are more complicated: we will show inductively that $\mathscr{C}_{k,i}$ is a $T_k \cdots T_i$ compact extension of $\mathscr{C}_{k-1,i}$ for each $1 \le i \le k-1$, and that $\mathscr{C}_{k,k}$ is a $T_k$ almost-periodic extension of $\mathscr{C}_{k-1} = \mathscr{C}_{k-1,1} \vee \cdots \vee \mathscr{C}_{k-1,k-1}$.
(See Figure \ref{fig:CharacteristicFactorsSchematic} on Page \pageref{fig:CharacteristicFactorsSchematic} for a schematic.)
It is not clear whether these characteristic factors can be used to prove Bergelson's conjecture.
The difficulty lies partly in their dependence on $i$ which, as exemplified in \cite{MR1357755}, cannot be removed in general.
Under the above-mentioned ergodicity assumptions we can handle this dependence when $k = 3$ and obtain multiple recurrence.

The rest of the paper runs as follows.
In the next two sections we recall definitions and results used throughout the remainder of the paper.
In Section \ref{sec:AlmostPeriodicFunctions} we prove some facts about almost-periodic functions and eigenfunctions over a factor that we will need to prove our factors are characteristic.
Section \ref{sec:characteristicFactors} contains a definition of the factors $\mathscr{C}_{k,i}$ and a proof that they are characteristic.
The following section contains a result that allows us to lift multiple recurrence from a single $\sigma$-algebra to a family of $\sigma$-algebras.
It is used in Section \ref{sec:recurrenceResults} to prove our multiple recurrence result.
Finally, we present some further consequences of our description of characteristic factors, including some combinatorial results, in Section \ref{sec:furtherResults}.

Thanks are due to the author's advisor, Vitaly Bergelson, for bringing the question to the author's attention and for his participation in many fruitful discussions while the paper was in preparation, to Angelo Nasca and Younghwan Son for their useful comments on the manuscript, and to Alexander Leibman for finding the reference \cite{MR1357755}.
We would also like to thank the anonymous referee for a patient, detailed report and for suggesting a streamlined proof of Theorem~\ref{thm:liftingCharacteristicFactors}.

\section{Preliminaries}
\label{sec:Preliminaries}

In this section we recall the facts we will need about measurable group actions, factors, disintegration of measures, joinings and $\ip^*$ sets.
We also give suitable versions of the van der Corput trick and the mean ergodic theorem.
For more details, see \cite{MR0498471}, \cite{MR603625} and \cite{MR1958753}.

Fix throughout this paper a locally-compact, second-countable, amenable group $G$ with a left Haar measure $\haar$ and a countable, dense subgroup $\Gamma$.
Amenability implies (4.16 in \cite{MR961261}) the existence of a sequence $\Phi$ of compact, positive-measure subsets of $G$ such that
\begin{equation*}
\frac{\haar(\Phi_N \symdiff g\Phi_N)}{\haar(\Phi_N)} \to 0
\end{equation*}
as $N \to \infty$ for each $g \in G$. 
The convergence is uniform on compact subsets of $G$.
Any such sequence is called a \define{(left) F{\o}lner sequence}.
Fix a left F{\o}lner sequence $\Phi$ in $G$.

Let $(X,\mathscr{B},\mu)$ be a separated, countably generated probability space.
By a \define{measurable action} of $G$ on such a space we mean a family $\{ T^g \,:\, g \in G \}$ of measurable, measure-preserving transformations of $(X,\mathscr{B},\mu)$ such that the induced map $G \times X \to X$ given by $(g,x) \mapsto T^g x$ is measurable and $T^gT^h = T^{gh}$ for all $g, h$ in $G$.
Two such actions $T_1$ and $T_2$ are said to \define{commute} if $T_1^g T_2^h = T_2^h T_1^g$ for all $g,h$ in $G$, and if they do $T_1^gT_2^g$ is also a measurable action of $G$ on $(X,\mathscr{B},\mu)$.

By a \define{system} we mean a tuple $(X,\mathscr{B},\mu,T)$ consisting of a measurable action $T$ of $G$ on a separated, countably generated probability space $(X,\mathscr{B},\mu)$.
We often write $\mathbf{X}$ for $(X,\mathscr{B},\mu,T)$ and $\lp^p(\mathbf{X})$ for the corresponding real space $\lp^p(X,\mathscr{B},\mu)$.
Given a system $\mathbf{X}$, each $T^g$ induces a unitary operator on $\lp^2(\mathbf{X})$ given by $(T^g f)(x) = f(T^g x)$.
It is immediate that $T^g(T^h f) = T^{hg} f$ for all $g,h$ in $G$.
Since $\mathscr{B}$ is countably generated the Hilbert space $\lp^2(\mathbf{X})$ is separable.
By 22.20(b) in \cite{MR551496} and the fact that $G \times X \to X$ is measurable, the map $g \mapsto T^g$ is continuous in the strong operator topology.

Given a sub-$\sigma$-algebra $\mathscr{C}$ of $\mathscr{B}$ and $f$ in $\lp^2(X,\mathscr{B},\mu)$ the \define{conditional expectation} of $f$ on $\mathscr{C}$, denoted $\condex{f}{\mathscr{C}}$, is the orthogonal projection of $f$ onto the closed subspace $\lp^2(X,\mathscr{C},\mu)$.
We say that a sub-$\sigma$-algebra $\mathscr{C}$ of $\mathscr{B}$ is \define{$T$ invariant} if $(T^g)^{-1}C \in \mathscr{C}$  for all $C \in \mathscr{C}$ and all $g \in G$.
When this is the case each $T^g$ commutes with the conditional expectation $\condex{\,\cdot\,}{\mathscr{C}}$.

We say that a system $\mathbf{Y} = (Y,\mathscr{D},\lambda,S)$ is a \define{factor} of $\mathbf{X} = (X,\mathscr{B},\mu,T)$, or that $\mathbf{X}$ is an \define{extension} of $\mathbf{Y}$, if there is a measurable, measure-preserving map $\pi : X \to Y$, called the \define{factor map}, that intertwines the actions $T$ and $S$, meaning that $\pi(T^g x) = S^g(\pi x)$ for all $x$ in $X$ and all $g$ in $G$.
We will usually abuse notation by writing $\mu$ for $\lambda$ and $T$ for $S$.
To any factor $\mathbf{Y}$ of $\mathbf{X}$ we can associate the $T$-invariant sub-$\sigma$-algebra $\pi^{-1}\mathscr{D}$ of $\mathscr{B}$.
We can use $\pi$ to identify $\lp^2(\mathbf{Y})$ with $\lp^2(X,\pi^{-1}\mathscr{D},\mu,T)$ isometrically.
This lets us think of $\condex{f}{\pi^{-1}\mathscr{D}}$ as an element of $\lp^2(\mathbf{Y})$, which we will denote $\condex{f}{\mathbf{Y}}$.

By Lemma 3.1 in \cite{MR1191743} any closed subspace of $\lp^2(\mathbf{X})$ that is a lattice and contains the constants is of the form $\lp^2(X,\mathscr{C},\mu)$ for some sub-$\sigma$-algebra $\mathscr{C}$ of $\mathscr{B}$.
If the subspace is $T$-invariant then so is $\mathscr{C}$.
Proposition 2.1 in \cite{MR0409770} lets us associate with any $T$-invariant sub-$\sigma$-algebra $\mathscr{C}$ of $\mathscr{B}$ a system $\mathbf{Y}$ and a $T$-invariant, full-measure set $X'$ in $\mathscr{B}$ such that $(X',\mathscr{B},\mu,T)$ is an extension of $\mathbf{Y}$ via a factor map $\pi : X' \to Y$.
Since the probability space defined by $X'$ is also separated and countably generated, we will not distinguish between $X'$ and $X$ hereafter.

A factor map $\pi : \mathbf{X} \to \mathbf{Y}$ gives rise to a \define{disintegration} of $\mu$ over $\mathbf{Y}$, which is a $\lambda$ almost-surely defined family $\{ \mu_y \,:\, y \in Y \}$ of probability measures on $(X,\mathscr{B})$ with the following properties.
\begin{enumerate}
\item  For any $\mathscr{B}$-measurable function $f$ that is square-integrable the map
\begin{equation*}
y \mapsto \int f \intd \mu_y
\end{equation*}
is defined $\lambda$ almost-surely and $\mathscr{D}$-measurable.
\item  For any $\mathscr{B}$-measurable function $f$ that is square-integrable
\begin{equation*}
\condex{f}{\mathbf{Y}}(y) = \int f \intd\mu_y
\end{equation*}
$\lambda$ almost-surely.
\item The group $G$ permutes the family $\mu_y$ in the sense that, for any $g \in G$ and any $\mathscr{B}$-measurable function $f$ that is square-integrable one has
\begin{equation*}
\int T^g f \intd\mu_y = \int f \intd\mu_{S^g y}
\end{equation*}
$\lambda$ almost-surely.
\end{enumerate}
Care is taken to speak of a function $f : X \to \mathbb{R}$ rather than an equivalence class of functions in $\lp^2(X,\mathscr{B},\mu)$ because the measures $\mu_y$ may well be singular with respect to $\mu$.
Although each integrable function $f : X \to \mathbb{R}$ defines an equivalence class in the space $\lp^2(X,\mathscr{B},\mu_y)$ for almost every $y$, changing $f$ on a set of $\mu$ measure zero may not preserve all of these classes.
Write $\langle\cdot,\cdot\rangle_y$ for the inner product on $\lp^2(X,\mathscr{B},\mu_y)$ and $\nbar \cdot \nbar_y$ for the corresponding norm.
Given an invariant sub-$\sigma$-algebra $\mathscr{D}$, by the disintegration of $\mu$ over $\mathscr{D}$ we mean the family of measure $\nu_x = \mu_{\pi x}$ where $\mu_y$ is an almost-surely defined disintegration of $\mu$ over a factor corresponding to $\mathscr{D}$.
By an abuse of notation we will write $\mu_x$ for $\nu_x$.

We now recall some basic facts about joinings.
Let $(X_i,\mathscr{B}_i,\mu_i), 1 \le i \le k$ be probability spaces and let $\pi_i$ be the projection from $X_1 \times \cdots \times X_k$ to $X_i$.
We say that a probability measure $\nu$ on $(X_1 \times \cdots \times X_k, \mathscr{B}_1 \otimes \cdots \otimes \mathscr{B}_k)$ is a \define{standard measure} if $\nu(\pi_i^{-1}B) = \mu_i(B)$ for all $B \in \mathscr{B}_i$ and all $1 \le i \le k$.
A sequence $\nu_n$ of standard measures is said to converge to a standard measure $\nu$ if
\begin{equation*}
\nu_n(B_1 \times \cdots \times B_k) \to \nu(B_1 \times \cdots \times B_k)
\end{equation*}
for all $B_i$ in $\mathscr{B}_i$.
A \define{joining} of systems $\mathbf{X}_1,\dots,\mathbf{X}_k$ is any system $\mathbf{X} = (X,\mathscr{B},\nu,T)$ where $X = X_1 \times \cdots \times X_k$, $\mathscr{B} = \mathscr{B}_1 \otimes \cdots \otimes \mathscr{B}_k$, $T^g = T_1^g \times \cdots \times T_k^g$ and $\nu$ is a standard measure that is $T$-invariant.
Given a factor $\mathbf{Y}_i$ of $\mathbf{X}_i$ for each $i$, we can consider the system $\mathbf{Y} = (Y,\mathscr{D},\eta,T)$ made from $\mathbf{X}$ by projecting $\nu$ onto the product $(Y,\mathscr{D})$ of the underlying measurable spaces $(Y_i,\mathscr{D}_i)$.
Call a joining $\mathbf{X}$ of the $\mathbf{X}_i$ a \define{conditional product joining} relative to the factors $\mathbf{Y}_i$ if 
\begin{equation}
\label{eqn:conditionalProductMeasure}
\int f_1 \otimes \cdots \otimes f_k \intd\nu = \int \condex{f_1}{\mathbf{Y}_1} \otimes \cdots \otimes \condex{f_k}{\mathbf{Y}_k} \intd\eta
\end{equation}
for all $f_i$ in $\lp^\infty(\mathbf{X}_i)$.
Here $f_1 \otimes \cdots \otimes f_k$ denotes the function mapping $(x_1,\dots,x_k)$ to $f_1(x_1)\cdots f_k(x_k)$.
We can re-write \eqref{eqn:conditionalProductMeasure} as
\begin{equation}
\label{eqn:conditionalProductMeasureDisintegration}
\nu = \int \mu_{1,y_1} \otimes \cdots \otimes \mu_{k,y_k} \intd\eta(y_1,\dots,y_k)
\end{equation}
if $\mu_{i,y_i}$ is the almost-surely defined disintegration of $\mu_i$ over $\mathbf{Y}_i$.

Let $T_1,\dots,T_k$ be commuting, measurable actions of $G$ on a probability space $(X,\mathscr{B},\mu)$.
Define a measure $\nu_k$ on $(X^{k+1},\mathscr{B}^{k+1})$ by
\begin{equation}
\label{FurstenbergJoining}
\int f_1 \otimes \cdots \otimes f_{k+1} \intd\nu_k = \lim_{N \to \infty} \frac{1}{\haar(\Phi_N)} \int\limits_{\Phi_N} \! \int f_{k+1} \cdot \prod_{i=1}^k T_k^g \cdots T_i^g f_i \intd\mu \intd\haar(g)
\end{equation}
for any $f_1,\dots,f_{k+1}$ in $\lp^\infty(X,\mathscr{B},\mu)$.
The existence of the limit is justified by Theorem 1.3 in \cite{arXiv:1111.7292}.
Using the fact that $\Phi$ is a F\o{}lner sequence, one can show that $\nu_k$ is
\begin{equation*}
T_k T_{k-1} \cdots T_1 \times \cdots \times T_k T_{k-1} \times T_k \times I
\end{equation*}
invariant.
Thus the measure $\nu_k$ yields a joining of the systems
\begin{equation*}
(X,\mathscr{B},\mu,T_k \cdots T_1),\dots,(X,\mathscr{B},\mu,T_k),(X,\mathscr{B},\mu,I)
\end{equation*}
called the \define{Furstenberg joining} of the actions $T_1,\dots,T_k$.

Given two systems $\mathbf{X}_1 = (X_1,\mathscr{B}_1,\mu_1,T_1)$ and $\mathbf{X}_2 = (X_2,\mathscr{B}_2,\mu_2,T_2)$ having a common factor $\mathbf{Y} = (Y,\mathscr{D},\mu,T)$ via factor maps $\pi_1$ and $\pi_2$ respectively, we can form their \define{relatively independent joining}
\begin{equation*}
\mathbf{X}_1 \times_\mathbf{Y} \mathbf{X}_2 = (X_1 \times X_2, \mathscr{B}_1 \otimes \mathscr{B}_2, \nu, T_1 \times T_2)
\end{equation*}
where $\nu$ is the measure defined by
\begin{equation*}
\int f_1 \otimes f_2 \intd\nu = \int \condex{f_1}{\mathbf{Y}} \cdot \condex{f_2}{\mathbf{Y}} \intd\mu
\end{equation*}
for all $f_1$ in $\lp^\infty(\mathbf{X}_1)$ and all $f_2$ in $\lp^\infty(\mathbf{X}_2)$.
The measure is supported on the set
\begin{equation*}
\{ (x_1,x_2) \,:\, \pi_1x_1 = \pi_2 x_2 \}
\end{equation*}
so $\mathbf{Y}$ is a factor of $\mathbf{X}_1 \times_\mathbf{Y} \mathbf{X}_2$ in an unambiguous way.
If $\mu_{1,y}$ and $\mu_{2,y}$ are the almost-surely defined disintegrations of $\mu_1$ and $\mu_2$ over $\mathbf{Y}$ then
\begin{equation}
\label{relIndepJoiningBase}
\nu = \int \mu_{1,y} \otimes \mu_{2,y} \intd\mu(y)
\end{equation}
is the disintegration of $\nu$ over $\mathbf{Y}$.
We also recall that
\begin{equation}
\label{eqn:relIndepJoiningSecondFactor}
\nu = \int \mu_{1,\pi_2x_2} \otimes \delta_{x_2} \intd\mu_2 (x_2)
\end{equation}
is the disintegration of $\nu$ over $\mathbf{X}_2$.

We will need some basic facts about $\ip$ sets.
Given a sequence $\phi$ in $G$ define
\begin{equation*}
\fp(\phi) = \left\{ \phi(i_1) \cdots \phi(i_k) \,:\, k \in \mathbb{N}, i_1 < \cdots < i_k \in \mathbb{N} \right\}
\end{equation*}
and call a subset of $G$ an \define{IP set} if it contains $\fp(\phi)$ for some sequence $\phi$ in $G$.
By Hindman's theorem (see Lemma~2.1 in \cite{MR1264030}) the property of being an IP set is partition regular.
A subset of $G$ is said to be $\ip^*$ if its intersection with every $\ip$ set is non-empty.
It follows from partition regularity (see Lemma~9.5 in \cite{MR603625}) that the intersection of two $\ip^*$ sets is also $\ip^*$.
Finally, note that every $\ip^*$ subset of $G$ has the property that finitely many of its left-shifts cover $G$.
This is because the complement of a set failing to have this property contains a right-shift of any finite set and therefore contains an $\ip$ set.
Thus every measurable $\ip^*$ set has positive lower density with respect to $\Phi$, where
\begin{equation*}
\lowerdens_\Phi(E) = \liminf_{N \to \infty} \frac{\haar(E \cap \Phi_N)}{\haar(\Phi_N)}
\end{equation*}
is the \define{lower density} of a measurable subset $E$ of $G$ with respect to $\Phi$.
Replacing $\liminf$ with $\limsup$ gives the \define{upper density} of $E$, denoted $\upperdens(E)$, and when $\upperdens(E) = \lowerdens(E)$ their common value, the \define{density} of $E$, is denoted $\dens(E)$.

We conclude this section with versions of the van der Corput trick and the mean ergodic theorem suitable for our needs.
The Hilbert space valued integrals below are always taken in the sense of Bochner.

\begin{proposition}[van der Corput trick]
\label{vanDerCorputTrick}
Let $\mathscr{H}$ be a separable Hilbert space and let $u : G \to \mathscr{H}$ be weakly measurable and uniformly bounded in norm. If
\begin{equation*}
\limsup_{H \to \infty} \frac{1}{\haar(\Phi_H)^2} \int\limits_{\Phi_H} \! \int\limits_{\Phi_H} \limsup_{N \to \infty} \left| \frac{1}{\haar(\Phi_N)} \int\limits_{\Phi_N} \langle u(hg), u(lg) \rangle \intd\haar(g) \right| \intd\haar(h) \intd\haar(l) = 0
\end{equation*}
then
\begin{equation}
\label{vdcStarter}
\bnbar \frac{1}{\haar(\Phi_N)} \int\limits_{\Phi_N} u(g) \intd\haar(g) \bnbar
\end{equation}
converges to 0 as $N \to \infty$.
\begin{proof}
Fix $\epsilon > 0$. First note that given any $H$ in $\mathbb{N}$ one has
\begin{equation*}
\lim_{N \to \infty} \bnbar \frac{1}{\haar(\Phi_N)} \int\limits_{\Phi_N} u(g) \intd\haar(g) - \frac{1}{\haar(\Phi_N)} \int\limits_{\Phi_N}\frac{1}{\haar(\Phi_H)} \int\limits_{\Phi_H} u(hg) \intd\haar(h) \intd\haar(g) \bnbar = 0
\end{equation*}
by the dominated convergence theorem and the fact that $\Phi$ is a left F{\o}lner sequence. By the Cauchy-Schwarz inequality
\begin{align*}
 & \bnbar \frac{1}{\haar(\Phi_N)} \int\limits_{\Phi_N}\frac{1}{\haar(\Phi_H)} \int\limits_{\Phi_H} u(hg) \intd\haar(h) \intd\haar(g) \bnbar^2\displaybreak[2]\\
 \le \, & \frac{1}{\haar(\Phi_N)} \int\limits_{\Phi_N} \bnbar \frac{1}{\haar(\Phi_H)} \int\limits_{\Phi_H} u(hg) \intd\haar(h) \bnbar^2 \intd\haar(g)\displaybreak[2]\\
 = \, & \frac{1}{\haar(\Phi_H)^2} \int\limits_{\Phi_H} \! \int\limits_{\Phi_H} \frac{1}{\haar(\Phi_N)} \int\limits_{\Phi_N} \langle u(hg), u(kg) \rangle \intd\haar(g) \intd\haar(h) \intd\haar(k)
\end{align*}
which allows us to relate \eqref{vdcStarter} to the hypothesis and obtain the desired result.
\end{proof}
\end{proposition}

\begin{proposition}[Mean ergodic theorem]
\label{meanErgodic}
Let $T$ be a measurable action of $G$ on a separated, countably generated probability space $(X,\mathscr{B},\mu)$. Let $\mathscr{I}$ be the sub-$\sigma$-algebra of $T$-invariant sets. Then
\begin{equation*}
\lim_{N \to \infty} \frac{1}{\haar(\Phi_N)} \int\limits_{\Phi_N} T^g f \intd\haar(g) = \condex{f}{\mathscr{I}}
\end{equation*}
in norm for all $f$ in $\lp^2(X,\mathscr{B},\mu)$.
\end{proposition}

This is Theorem~5.7 in \cite{MR961261}.
In particular we have
\begin{equation*}
\lim_{N \to \infty} \frac{1}{\haar(\Phi_N)} \int\limits_{\Phi_N} \! \int T^g f_1 \cdot f_2 \intd\mu \intd\haar(g) = \int \condex{f_1}{\mathscr{I}} \cdot \condex{f_2}{\mathscr{I}} \intd\mu
\end{equation*}
for all $f_1,f_2$ in $\lp^2(X,\mathscr{B},\mu)$.

\section{Borel Hilbert bundles}
\label{sec:BorelHilbertBundles}

In this section we recall how to associate Borel Hilbert bundles with extensions and relatively independent joinings. For details, see \cite{MR641217}, \cite{MR1958753} and \cite{MR2288954}.

Let $(Y,\mathscr{D},\mu)$ be a separable, countably generated probability space and let $\mathfrak{H} = \{ \mathfrak{H}_y \,:\, y \in Y \}$ be a collection of separable, real Hilbert spaces.
Write $\langle \cdot,\cdot \rangle_y$ for the inner product on $\mathfrak{H}_y$.
From $Y$ and $\mathfrak{H}$ we can form the total space $Y \total \mathfrak{H} = \{ (y,h) \,:\, y \in Y, h \in \mathfrak{H}_y \}$ which comes with a projection $\pi : Y \total \mathfrak{H} \to Y$.
The spaces $\mathfrak{H}_y$ are called the \define{fibers} of the total space.
A \define{section} of $Y \total \mathfrak{H}$ is any map $f:Y \to Y \total \mathfrak{H}$ such that $\pi \circ f$ is the identity.
The image of a point $y$ under a section $f$ is a point in $Y \total \mathfrak{H}$ which we will write as $(y,f_y)$.
Thus $f_y$ belongs to $\mathfrak{H}_y$.
To any section $f$ we can associate the map $\tilde{f} : Y \total \mathfrak{H} \to \mathbb{R}$ defined by $\tilde{f}(y,h) = \langle f_y, h \rangle_y$.
A \define{Borel Hilbert bundle} is a Hilbert bundle $Y \total \mathfrak{H}$ equipped with a $\sigma$-algebra of subsets of $Y \total \mathfrak{H}$ for which:
\begin{enumerate}
\item[(i)] the projection $Y \total \mathfrak{H} \to Y$ is measurable;
\item[(ii)] there is a sequence $f_n$ of sections such that:\begin{enumerate}
\item the maps $\tilde{f_n}$ are measurable;
\item for each $n,m$ the map $Y \to \mathbb{R}$ given by $y \mapsto \langle f_{n,y}, f_{m,y} \rangle_y$ is measurable;
\item the functions $\tilde{f}_n$ and $\pi$ separate points on $Y \total \mathfrak{H}$.
\end{enumerate}
\end{enumerate}

To associate a Borel Hilbert bundle $Y \total \mathfrak{H}$ with a given extension $\mathbf{X} \to \mathbf{Y}$, fix an almost-surely defined disintegration $\mu_y$ of $\mu$ over $\mathbf{Y}$ and let $\mathscr{A} = \{ A_1, A_2,\dots \}$ be a countable, $\Gamma$-invariant sub-algebra of $\mathscr{B}$ that generates $\mathscr{B}$.
For each $n,m$ the function $y \mapsto \langle 1_{A_n}, 1_{A_m} \rangle_y$ is defined on a full-measure subset of $Y$ and is measurable there.
Let $Y_0$ be a $\Gamma$-invariant, full-measure subset of $Y$ on which $\mu_y$ and all of the functions $y \mapsto \langle 1_{A_n}, 1_{A_m} \rangle_y$ are defined and on which $T^\gamma \mu_y = \mu_{T^\gamma y}$ for all $\gamma$ in $\Gamma$.
Put $\mathfrak{H}_y = \lp^2(X,\mathscr{B},\mu_y)$ when $y \in Y_0$ and put $\mathfrak{H}_y = \{ 0 \}$ otherwise.
Each $\mathfrak{H}_y$ is separable because $\mathscr{B}$ is countably generated.
Let $\mathfrak{H}$ be the collection $\{ \mathfrak{H}_y \,:\, y \in Y \}$.
Define a sequence $f_n$ of sections by taking $f_{n,y} = 1_{A_n}$ when $y \in Y_0$ and $f_{n,y} = 0$ otherwise.
Equip $Y \total \mathfrak{H}$ with the smallest $\sigma$-algebra of subsets for which $\pi$ and the maps $\tilde{f}_n$ are measurable.
It is immediate from the construction that this $\sigma$-algebra makes $Y \total \mathfrak{H}$ into a Borel Hilbert bundle.
Moreover, a section $f : Y \to Y \total \mathfrak{H}$ is measurable with respect to this $\sigma$-algebra if and only if $y \mapsto \langle f_y, f_{n,y} \rangle_y$ is measurable for each $n$.
We call $Y \total \mathfrak{H}$ the Borel Hilbert bundle corresponding to the extension $\mathbf{X} \to \mathbf{Y}$.
The Hilbert space $\lp^2(Y \total \mathfrak{H},\mu)$ formed from the set
\begin{equation*}
\mathscr{L}^2(Y \total \mathfrak{H},\mu) = \{ f \in B(Y \total \mathfrak{H}) \,:\, y \mapsto \nbar f_y \nbar^2_y \text{ is } \mu \text{ integrable} \}
\end{equation*}
of square-integrable sections by identifying sections that agree almost surely is isomorphic to $\lp^2(X,\mathscr{B},\mu)$.
Thus to any $\phi$ in $\lp^2(\mathbf{X})$ we can associate an almost-surely defined, square-integrable section $y \mapsto \phi_y$ and vice versa.

We now recall how $\Gamma$ acts on sections of $Y \total \mathfrak{H}$.
Fix $\gamma \in \Gamma$.
Since $T^\gamma \mu_y = \mu_{T^\gamma y}$ whenever $y \in Y_0$ the map $T_y^\gamma : \mathfrak{H}_{T^\gamma y} \to \mathfrak{H}_y$ given by $(T_y^\gamma f)(x) = f(T^\gamma x)$ is well-defined and unitary.
Define $T_y^\gamma : \mathfrak{H}_{T^\gamma y} \to \mathfrak{H}_y$ to be the zero map when $y \notin Y_0$.
The family of maps $\{ T^\gamma_y \,:\, y \in Y \}$ induces a map $T^\gamma$ on sections of $Y \total \mathfrak{H}$ such that $(T^\gamma f)_y = T^\gamma_y f_{T^\gamma y}$.
If $f$ is a measurable section then so is $T^\gamma f$.
Also $T^\eta(T^\gamma f) = T^{\gamma\eta} f$ for all $\gamma, \eta$ in $\Gamma$.

It remains to relate the Borel Hilbert bundle associated with a relatively independent joining $\mathbf{X}_1 \times_\mathbf{Y} \mathbf{X}_2 \to \mathbf{Y}$ to the Hilbert bundles associated with the extensions $\mathbf{X}_1 \to \mathbf{Y}$ and $\mathbf{X}_2 \to \mathbf{Y}$.
Let $\mathscr{A}_1$ and $\mathscr{A}_2$ be countable, $\Gamma$-invariant algebras that generate $\mathscr{B}_1$ and $\mathscr{B}_2$ respectively.
The countable algebra generated by $\mathscr{A}_1 \otimes \mathscr{A}_2$ is $\Gamma$-invariant and generates $\mathscr{B}_1 \otimes \mathscr{B}_2$.
We can thus simultaneously form the Borel Hilbert bundles $Y \total \mathfrak{H}_1$, $Y \total \mathfrak{H}_2$ and $Y \total \mathfrak{H}$ corresponding to the extensions $\mathbf{X}_1 \to \mathbf{Y}$, $\mathbf{X}_2 \to \mathbf{Y}$ and $\mathbf{X}_1 \times_\mathbf{Y} \mathbf{X}_2 \to \mathbf{Y}$ respectively.
From \eqref{relIndepJoiningBase} we see that $\mathfrak{H}_y = \mathfrak{H}_{1,y} \otimes \mathfrak{H}_{2,y}$ for $\mu$ almost every $y$.
Thus for any section $H$ of $Y \total \mathfrak{H}$ and almost every $y$ the corresponding member $H_y$ of $\mathfrak{H}_y$ induces a compact operator $H_y : \mathfrak{H}_{1,y} \to \mathfrak{H}_{2,y}$ defined by
\begin{equation}
\label{bundleConvolution}
(\phi \conv H_y)(x_2) = \int \phi_1(x_1) \cdot H(x_1,x_2) \intd\mu_{1,y}(x_1)
\end{equation}
for any $\phi \in \mathfrak{H}_{1,y}$.
This family of operators induces a map taking almost-surely defined sections of $Y \total \mathfrak{H}_1$ to almost-surely defined sections of $Y \total \mathfrak{H}_2$.
If $H$ is a measurable section of $Y \total \mathfrak{H}$ the induced map preserves measurability of sections because
\begin{equation*}
y \mapsto \langle 1_{A_{1,n}} \conv H, 1_{A_{2,m}} \rangle_y = \langle H, 1_{A_{1,n}} \otimes 1_{A_{2,m}} \rangle_y
\end{equation*}
is measurable for all $n,m$ in $\mathbb{N}$.
However, the induced map need not preserve square-integrability.
It may happen that $f_1$ is a square-integrable section of $Y \total \mathfrak{H}_1$ and $y \mapsto f_{1,y} \conv H_y$ is not a square-integrable section of $Y \total \mathfrak{H}_2$.
As the following proposition shows, we avoid this problem when the norms $\nbar H_y \nbar_y$ are bounded almost-surely and write $f_1 \conv H$ for the element of $\lp^2(\mathbf{X}_2)$ corresponding to the square integrable section $f_{1,y} \conv H_y$ of $Y \total \mathfrak{H}_2$.

\begin{proposition}
\label{prop:sectionOfBoundedOperators}
Let $H$ be a section of $Y \total \mathfrak{H}$.
If the norms $\nbar H_y \nbar_y$ are essentially bounded and $f_1$ is a square-integrable section of the bundle $Y \total \mathfrak{H}_1$ then $f_{1,y} \conv H_y$ is a square-integrable section of $Y \total \mathfrak{H}_2$.
\begin{proof}
See Section F.3 in \cite{MR2288954}.
\end{proof}
\end{proposition}

A section $H$ of $Y \total \mathfrak{H}$ also defines for almost every $y$ a compact operator $H_y : \mathfrak{H}_{2,y} \to \mathfrak{H}_{1,y}$ defined by
\begin{equation*}
(H_y \conv \phi)(x_1) = \int H_y(x_1,x_2) \cdot \phi(x_2) \intd\mu_{2,y}(x_2)
\end{equation*}
for any $\phi \in \mathfrak{H}_{2,y}$ with similar properties.

Given a measurable section $H$ of $Y \total \mathfrak{H}$ we can spectrally decompose the compact operator $H_y : \mathfrak{H}_{1,y} \to \mathfrak{H}_{2,y}$ for almost every $y \in Y$. The following theorem, due to Furstenberg and Katznelson, shows that when $\mathbf{X}_1 = \mathbf{X}_2$ and $H$ is positive-definite and symmetric, the spectral decomposition is measurable.

\begin{theorem}[3.7 in \cite{MR1191743}]
\label{joiningEigensections}
Let $\mathbf{X} \to \mathbf{Y}$ be an extension of systems. Form the corresponding Borel Hilbert bundle $Y \total \mathfrak{H}$. Let $H_y$ be a measurable family of positive-definite, self-adjoint, compact operators on $\mathfrak{H}_y$. Let $\lambda_n(y)$ be a decreasing enumeration of the positive eigenvalues of $H_y$, counting multiplicities. There is a sequence $\Psi_n$ of square integrable sections of $Y \total \mathfrak{H}$ such that $\Psi_{n,y} \conv H_y = \lambda_n(y)\Psi_{n,y}$ whenever $\lambda_n(y)$ is defined, $\Psi_{n,y} = 0$ otherwise, and $\{ \Psi_{n,y} \,:\, n \in \mathbb{N} \} \backslash \{ 0 \}$ is orthonormal in almost every fiber.
\end{theorem}

\section{Almost-Periodic Functions and Eigenfunctions}
\label{sec:AlmostPeriodicFunctions}

We will describe the characteristic factors $\mathscr{C}_{k,i}$ in terms of almost periodic functions.
In this section we prove the results about almost-periodic functions and eigenfunctions that we will need later.
Most of the results in this section are well-known in one form or another; we provide the details for the sake of completion.

Let $\mathbf{X} \to \mathbf{Y}$ be an extension and let $\mu_y$ be an almost-surely defined disintegration of $\mu$ over $\mathbf{Y}$.
We say that $f$ in $\lp^2(\mathbf{X})$ is \define{almost-periodic} for this extension if for every $\epsilon > 0$ one can find a finite subset $\Xi$ of $\lp^\infty(\mathbf{X})$ and $E \subset Y$ with $\mu(E) > 1 - \epsilon$ such that
\begin{equation}
\label{def:AlmostPeriodic}
\min \{ \nbar T^\gamma f - \xi \nbar_y \,:\, \xi \in \Xi \} \le \epsilon
\end{equation}
for each $\gamma \in \Gamma$ and almost every $y \in E$.
The closure of the set of almost-periodic functions, which we denote $\ap(\mathbf{X}|\mathbf{Y})$, forms a closed subspace of $\lp^2(\mathbf{X})$ that contains the constant functions.
Also, if $f$ is almost-periodic then so is $|f|$.
Thus condition (c) in Lemma 3.1 of \cite{MR1191743} is satisfied and there exists a sub-$\sigma$-algebra $\mathscr{C}$ of $\mathscr{B}$ such that $\ap(\mathbf{X}|\mathbf{Y}) = \lp^2(X,\mathscr{C},\mu)$.
This lets us approximate any $f$ in $\ap(\mathbf{X}|\mathbf{Y})$ arbitrarily well by a function in $\lp^\infty(X,\mathscr{B},\mu)$ that is almost-periodic over $\mathbf{Y}$ as follows: truncate $f$ at a high level and then re-define $f$ to be zero on certain fibers of the factor map as in the proof of Theorem~9.1 in \cite{MR670131}.
Since $\mathcal{A}(\mathbf{X}|\mathbf{Y})$ is closed and invariant under $T^\gamma$ for each $\gamma \in \Gamma$, it is also $T$ invariant.
Thus $\mathscr{C}$ is $T$ invariant.
When $\mathbf{Y}$ is the trivial factor, write $\ap(\mathbf{X})$ for $\ap(\mathbf{X}|\mathbf{Y})$.

We say that an extension $\mathbf{X} \to \mathbf{Y}$ is \define{compact} if $\ap(\mathbf{X}|\mathbf{Y}) = \lp^2(\mathbf{X})$ and \define{weak-mixing} if $\ap(\mathbf{X}|\mathbf{Y}) = \lp^2(\mathbf{Y})$.
Given sub-$\sigma$-algebras $\mathscr{D}$ and $\mathscr{E}$ of $\mathscr{B}$, we will say that $\mathscr{D} \to \mathscr{E}$ is \define{compact} for $T$ if $\mathscr{D}$ and $\mathscr{E}$ are $T$ invariant sub-$\sigma$-algebras of $\mathscr{B}$ and the corresponding extension $\mathbf{Y} \to \mathbf{Z}$ is compact.

Fix systems $\mathbf{X}_1$ and $\mathbf{X}_2$ having a common factor $\mathbf{Y}$.
Let $\mathscr{C}_1$ and $\mathscr{C}_2$ be the $\sigma$-algebras corresponding to $\ap(\mathbf{X}_1|\mathbf{Y})$ and $\ap(\mathbf{X}_2|\mathbf{Y})$ respectively.
We begin by relating $\ap(\mathbf{X}_1|\mathbf{Y})$ and $\ap(\mathbf{X}_2|\mathbf{Y})$ to $\ap(\mathbf{X}_1 \times_\mathbf{Y} \mathbf{X}_2|\mathbf{Y})$ by showing that any $H$ in $\lp^\infty(\mathbf{X}_1 \times_\mathbf{Y} \mathbf{X}_2)$ that is almost-periodic over $\mathbf{Y}$ satisfies
\begin{equation}
\label{relIndepJoining}
\langle H,f_1 \otimes f_2 \rangle = \langle H, \condex{f_1}{\mathscr{C}_1} \otimes \condex{f_2}{\mathscr{C}_2} \rangle
\end{equation}
for any $f_1$ in $\lp^\infty(\mathbf{X}_1)$ and any $f_2$ in $\lp^\infty(\mathbf{X}_2)$.
This is similar to Proposition 4.4.4 in \cite{MR1738544}.

\begin{proposition}
\label{convolveWithApGivesAp}
For any $H$ in $\lp^\infty(\mathbf{X}_1 \times_\mathbf{Y} \mathbf{X}_2)$ almost-periodic over $\mathbf{Y}$ and any $f_1$ in $\lp^2(\mathbf{X}_1)$ the element $f_1 \conv H$ of $\lp^2(\mathbf{X}_2)$ is almost-periodic over $\mathbf{Y}$.
\begin{proof}
It suffices to prove this when $f$ is in $\lp^\infty(\mathbf{X}_1)$.
Fix $\epsilon > 0$.
We have to find a finite subset $\Xi$ of $\lp^\infty(\mathbf{X}_2)$ and a subset $E$ of $Y$ with $\mu(E) > 1 - \epsilon$ such that
\begin{equation}
\label{apInvariantGoal}
\min \{ \nbar T_2^\gamma (f_1 \conv H) - \xi \nbar_y \,:\, \xi \in \Xi \} \le \epsilon
\end{equation}
for each $\gamma \in \Gamma$ and almost every $y \in E$.
Almost-periodicity of $H$ over $\mathbf{Y}$ implies the existence of a finite subset $\Psi$ of $\lp^\infty(\mathbf{X}_1 \times_\mathbf{Y} \mathbf{X}_2)$ and a subset $F_1$ of $Y$ with $\nu(F_1) > 1 - \epsilon/16$ such that
\begin{equation*}
\min \{ \nbar (T_1^\gamma \times T_2^\gamma) H - \psi \nbar_y \,:\, \psi \in \Psi \} < \epsilon/16
\end{equation*}
for all $\gamma \in \Gamma$ and almost all $y \in F_1$.
Write $\Psi = \{ \psi_1,\dots,\psi_k \}$. Let $\gamma_n$ be an enumeration of $\Gamma$.
For each $1 \le i \le k$ and almost every $y$ we have a compact operator $\psi_{i,y}$ mapping $\mathfrak{H}_{1,y}$ to $\mathfrak{H}_{2,y}$.
Thus for each $1 \le i \le k$ and almost every $y$ we can find a positive integer $M_i(y)$ such that
\begin{equation*}
\{ (T_1^{\gamma_n} f_{1,T^{\gamma_n}y}) \conv \psi_{i,y}\,:\, 1 \le n \le M_i(y) \}
\end{equation*}
is $\epsilon/16$-dense in $\{ (T_1^\gamma f_{1,T^\gamma y}) \conv \psi_{i,y} \,:\, \gamma \in \Gamma \}$.
Each of the functions $M_i$ is measurable.
Thus we can find $N$ in $\mathbb{N}$ so large that $F_2 = M_1^{-1}[1,N] \cap \cdots \cap M_k^{-1}[1,N]$ has measure at least $1 - \epsilon/16$. Put
\begin{equation*}
\Xi = \{ (T^{\gamma_n} f_1) \conv \psi_i \,:\, 1 \le i \le k, 1 \le n \le N \}
\end{equation*}
and $E = F_1 \cap F_2$.
Fix $\gamma$ in $\Gamma$ and almost any $y$ in $E$.
We can choose $i$ such that
\begin{equation*}
\nbar (T_1^\gamma \times T_2^\gamma)H - \psi_i\nbar_y \le \epsilon/16
\end{equation*}
and then guarantee
\begin{equation*}
\nbar (T_1^\gamma f_{1,T^\gamma y}) \conv \psi_{i,y} - (T_1^{\gamma_n} f_{1,T^{\gamma_n} y}) \conv \psi_{i,y} \nbar_y \le \epsilon/16
\end{equation*}
holds for some $1 \le n \le N$.
From
\begin{align*}
(T_2^\gamma(f_1 \conv H))_y(x_2) & = \int f_1(x_1) H(x_1,T_2^\gamma x_2) \intd\mu_{1,T^\gamma y}(x_1)\\
 & = \int f_1(T_1^\gamma x_1) H(T_1^\gamma x_1, T_2^\gamma x_2) \intd\mu_{1,y}(x_1)
\end{align*}
we have
\begin{align*}
 & \nbar T_2^\gamma(f_1 \conv H) - (T_1^{\gamma_n} f_1) \conv \psi_i \nbar_y\\
 \le\, & \nbar T_2^\gamma(f_1 \conv H) - (T_1^\gamma f_1) \conv  \psi_i \nbar_y + \nbar (T_1^\gamma f_1) \conv \psi_i - (T_1^{\gamma_n} f_1) \conv \psi_i \nbar_y\\
 \le\, & \nbar (T_1^\gamma f_1) \conv ((T_1^\gamma \times T_2^\gamma)H) - (T_1^\gamma f_1) \conv \psi_i \nbar_y + \nbar (T_1^\gamma f_1) \conv \psi_i - (T_1^{\gamma_n} f_1) \conv \psi_i \nbar_y
\end{align*}
so \eqref{apInvariantGoal} holds as desired.
\end{proof}
\end{proposition}

\begin{proposition}
For any $H$ in $\lp^\infty(\mathbf{X}_1 \times_\mathbf{Y} \mathbf{X}_2)$ almost-periodic over $\mathbf{Y}$ and any $f_2$ in $\lp^2(\mathbf{X}_2)$ the element $H \conv f_2$ of $\lp^2(\mathbf{X}_1)$ is almost-periodic over $\mathbf{Y}$.
\begin{proof}
Identical to the proof of the previous proposition.
\end{proof}
\end{proposition}

\begin{proposition}
\label{ApLivesOnAps}
For any $H$ in $\lp^\infty(\mathbf{X}_1 \times_\mathbf{Y} \mathbf{X}_2)$ that is almost-periodic over $\mathbf{Y}$, any $f_1$ in $\lp^2(\mathbf{X}_1)$ and any $f_2$ in $\lp^2(\mathbf{X}_2)$ we have \eqref{relIndepJoining}.
\begin{proof}
We have
\begin{equation*}
\langle H, f_1 \otimes f_2 \rangle = \langle H, (f_1 - \condex{f_1}{\mathscr{C}_1} + \condex{f_1}{\mathscr{C}_1}) \otimes f_2 \rangle 
\end{equation*}
and a similar equality holds for $f_2$ so it suffices to prove that $\langle H, f_1 \otimes f_2 \rangle$ is zero when either $f_1$ is orthogonal to $\mathscr{C}_1$ or $f_2$ is orthogonal to $\mathscr{C}_2$.
The two cases are similar. In the latter we have
\begin{align*}
\langle H, f_1 \otimes f_2 \rangle &= \iint f_1(x_1)H(x_1,x_2)f_2(x_2) \intd(\mu_{1,y} \otimes \mu_{2,y})(x_1,x_2) \intd\mu(y)\\
 &= \iint (f_1 \conv H)(x_2) f_2(x_2) \intd\mu_{2,y}(x_2) \intd\mu(y) = \langle f_1 \conv H, f_2 \rangle
\end{align*}
which is zero by Proposition \ref{convolveWithApGivesAp}.
\end{proof}
\end{proposition}

\begin{corollary}
\label{invariantWmIsZero}
Let $\mathscr{I}$ be the sub-$\sigma$-algebra of $T_1 \times T_2$ invariant sets in $\mathbf{X}_1 \times_\mathbf{Y} \mathbf{X}_2$.
If $f_1$ in $\lp^\infty(\mathbf{X}_1)$ is orthogonal to $\ap(\mathbf{X}_1|\mathbf{Y})$ or $f_2$ in $\lp^\infty(\mathbf{X}_2)$ is orthogonal to $\ap(\mathbf{X}_2|\mathbf{Y})$ then $\condex{f_1 \otimes f_2}{\mathscr{I}} = 0$ in $\lp^2(\mathbf{X}_1 \times_\mathbf{Y} \mathbf{X}_2)$.
\begin{proof}
For any $T_1 \times T_2$ invariant function $H$ in $\lp^2(\mathbf{X}_1 \times_\mathbf{Y} \mathbf{X}_2)$ we have
\begin{equation*}
\langle \condex{f_1 \otimes f_2}{\mathscr{I}}, H \rangle = \langle f_1 \otimes f_2, H \rangle = \langle \condex{f_1}{\mathscr{C}_1} \otimes \condex{f_2}{\mathscr{C}_2}, H \rangle = 0
\end{equation*}
by \eqref{relIndepJoining}, because invariant functions are certainly almost-periodic.
\end{proof}
\end{corollary}

\begin{theorem}
\label{apIndepJoiningOverBase}
$\ap(\mathbf{X}_1 \times_\mathbf{Y} \mathbf{X}_2|\mathbf{Y}) = \ap(\mathbf{X}_1|\mathbf{Y}) \otimes \ap(\mathbf{X}_2|\mathbf{Y})$.
\begin{proof}
It is straightforward to check that if $f_1$ in $\lp^\infty(\mathbf{X}_1)$ and $f_2$ in $\lp^\infty(\mathbf{X}_2)$ are both almost-periodic over $\mathbf{Y}$ then $f_1 \otimes f_2$ in $\lp^2(\mathbf{X}_1 \times_\mathbf{Y} \mathbf{X}_2)$ is almost-periodic over $\mathbf{Y}$.
On the other hand, if $H$ belongs to $\ap(\mathbf{X}_1\times_\mathbf{Y}\mathbf{X}_2|\mathbf{Y})$ and is orthogonal to $\ap(\mathbf{X}_1|\mathbf{Y}) \otimes \ap(\mathbf{X}_2|\mathbf{Y})$ then by Proposition \ref{ApLivesOnAps} we have $\langle H, f_1 \otimes f_2 \rangle = 0$ for all $f_1$ in $\lp^\infty(\mathbf{X}_1)$ and all $f_2$ in $\lp^\infty(\mathbf{X}_2)$ so $H = 0$.
\end{proof}
\end{theorem}

Recall that a function $f$ in $\lp^2(\mathbf{X})$ is \define{weakly mixing} for $\mathbf{X} \to \mathbf{Y}$ if
\begin{equation*}
\lim_{N \to \infty} \frac{1}{\haar(\Phi_N)} \int\limits_{\Phi_N} \! \int | \condex{\phi \cdot T^g f}{\mathbf{Y}} |^2 \intd\mu \intd\haar(g) = 0
\end{equation*}
for every $\phi$ in $\lp^\infty(\mathbf{X})$. The set $\mathcal{W}(\mathbf{X}|\mathbf{Y})$ of weakly mixing functions is a closed, $T$ invariant subspace of $\lp^2(\mathbf{X})$.
Proposition \ref{convolveWithApGivesAp} lets us prove the following result.

\begin{theorem}
\label{thm:ApWmOrthogonal}
For any extension $\mathbf{X} \to \mathbf{Y}$ we have $\lp^2(\mathbf{X}) = \ap(\mathbf{X}|\mathbf{Y}) \oplus \wm(\mathbf{X}|\mathbf{Y})$.
\begin{proof}
First we show that if $f$ in $\lp^2(\mathbf{X})$ is orthogonal to $\ap(\mathbf{X}|\mathbf{Y})$ then $f$ belongs to $\wm(\mathbf{X}|\mathbf{Y})$.
Fix $\phi$ in $\lp^\infty(\mathbf{X})$.
Let $\mathscr{I}$ denote the sub-$\sigma$-algebra of $T \times T$ invariant sets and put $H = \condex{\phi \otimes \phi}{\mathscr{I}}$ in $\lp^2(\mathbf{X} \times_\mathbf{Y} \mathbf{X})$.
We have
\begin{equation*}
\begin{aligned}
 & \lim_{N \to \infty} \frac{1}{\haar(\Phi_N)} \int\limits_{\Phi_N} \! \int | \condex{\phi \cdot T^g f}{\mathbf{Y}} |^2 \intd\mu \intd\haar(g)\\
= & \lim_{N \to \infty} \frac{1}{\haar(\Phi_N)} \int\limits_{\Phi_N} \! \int (\phi \otimes \phi) \cdot (T \times T)^g (f \otimes f) \intd\nu \intd\haar(g)\\
= & \int \condex{\phi \otimes \phi}{\mathscr{I}} \cdot (f \otimes f) \intd\nu = \langle H, f \otimes f\rangle = \langle f \conv H, f\rangle = 0
\end{aligned}
\end{equation*}
by the mean ergodic theorem and Proposition~\ref{convolveWithApGivesAp}.
Since $\phi$ was arbitrary, $f$ is weakly mixing over $\mathbf{Y}$.

Now we show that $\ap(\mathbf{X}|\mathbf{Y})$ and $\wm(\mathbf{X}|\mathbf{Y})$ are orthogonal.
Fix $f$ in $\wm(\mathbf{X}|\mathbf{Y})$.
It suffices to prove that $f$ is orthogonal to any $\phi$ in $\lp^\infty(\mathbf{X})$ that is almost-periodic over $\mathbf{Y}$. 
Fix $\epsilon > 0$.
Since $\phi$ is almost periodic we can find a subset $E$ of $Y$ with $\mu(E) > 1 - \epsilon$ and a finite subset $\Xi = \{ \xi_1,\dots,\xi_k \}$ of $\lp^2(\mathbf{X})$ such that \eqref{def:AlmostPeriodic} holds for all $\gamma \in \Gamma$ and all $y \in E$.
Fix $g \in G$.
Since $\Gamma$ is dense in $G$ we can find some $\gamma$ in $\Gamma$ such that $\nbar T^g \phi - T^\gamma \phi \nbar_y < \epsilon$ for all $y$ in a subset $E_g$ of $Y$ with $\mu(E_g) > 1 - \epsilon$.
For each $y$ in $E$ choose $1 \le \iota(y) \le k$ so that $\nbar T^\gamma \phi - \xi_{\iota(y)}\nbar_y < \epsilon$.
Put $F = E \cap E_g$.
Cauchy-Schwarz gives
\begin{equation}
\label{eqn:wmApOrtho1}
\begin{aligned}
\int T^g \phi \cdot T^g f \intd\mu_y & \le \left| \int \xi_{\iota(y)} \cdot T^g f \intd\mu_y \right| + 2\epsilon \nbar T^g f \nbar_y\\
 & \le \sum_{i=1}^k | \condex{\xi_i \cdot T^g f}{\mathbf{Y}}(y) | + 2 \epsilon \nbar T^g f \nbar_y
\end{aligned}
\end{equation}
for any $y \in F$. Combining this with
\begin{equation}
\label{eqn:wmApOrtho2}
\int T^g \phi \cdot T^g f \intd\mu_y \le \int |T^g f| \intd\mu_y \nbar \phi \nbar_\infty \le \nbar T^g f \nbar_y \nbar \phi \nbar_\infty
\end{equation}
which holds (in particular) for almost-every $y \notin F$ we get
\begin{equation*}
\begin{aligned}
|\langle \phi,f \rangle| & \le \sum_{i=1}^k \int |\condex{\xi_i \cdot T^g f}{\mathbf{Y}}| \intd\mu + 2 \epsilon \nbar T^g f \nbar + \int 1_{Y \backslash F}(y) \cdot \nbar T^g f \nbar_y \intd\mu(y) \nbar \phi \nbar_\infty\\
 & \le \sum_{i=1}^k \int |\condex{\xi_i \cdot T^g f}{\mathbf{Y}}| \intd\mu + 2\epsilon \nbar f \nbar + \sqrt{2\epsilon} \cdot \nbar f \nbar \cdot \nbar \phi \nbar_\infty
\end{aligned}
\end{equation*}
by integrating, applying Cauchy-Schwarz, and noting that $\mu(Y \backslash F) \le 2\epsilon$.
Finally, averaging over the F\o{}lner sequence $\Phi$ and applying Cauchy-Schwarz once more gives
\begin{equation*}
|\langle \phi, f \rangle| \le 2\epsilon \nbar f \nbar + \sqrt{2\epsilon} \nbar f \nbar \nbar \phi \nbar_\infty + \sum_{i=1}^k \bigg( \frac{1}{\haar(\Phi_N)} \int\limits_{\Phi_N} \! \int |\condex{\xi_i \cdot T^g f}{\mathbf{Y}}|^2 \intd\mu \intd\haar(g) \bigg)^{1/2}
\end{equation*}
which, upon using the fact that $f$ is weakly-mixing and noting that $\epsilon$ was arbitrary, gives $\langle \phi, f \rangle = 0$.
\end{proof}
\end{theorem}

Since the definition of $\ap(\mathbf{X}|\mathbf{Y})$ is independent of the F\o{}lner sequence $\Phi$, the above proposition implies that $\wm(\mathbf{X}|\mathbf{Y})$ is also independent of $\Phi$.

We will use Theorem~\ref{thm:ApWmOrthogonal} to relate $\ap(\mathbf{X}|\mathbf{Y})$ to the eigenfunctions of an extension.
Given an extension $\mathbf{X} \to \mathbf{Y}$, the factor map lets us embed $\lp^\infty(\mathbf{Y})$ in $\lp^\infty(\mathbf{X})$.
Thus we can think of $\lp^2(\mathbf{X})$ as an $\lp^\infty(\mathbf{Y})$ module.
A function $f$ in $\lp^2(\mathbf{X})$ is an \define{eigenfunction} over $\mathbf{Y}$ if the closed subspace $\mathscr{M}$ spanned by the orbit of $f$ is a finite-rank $\lp^\infty(\mathbf{Y})$ module.
This means we can find $\phi_1,\dots,\phi_d$ in $\lp^2(\mathbf{X})$ such that
\begin{equation*}
\{ \alpha^1 \phi_1 + \cdots + \alpha^d \phi_d \,:\, \alpha^1,\dots,\alpha^d \in \lp^\infty(\mathbf{Y}) \}
\end{equation*}
is dense in $\mathscr{M}$.
Denote by $\eig(\mathbf{X}|\mathbf{Y})$ the closed subspace of $\lp^2(\mathbf{X})$ spanned by the eigenfunctions over $\mathbf{Y}$.
When $\mathbf{Y}$ is the trivial factor, write $\eig(\mathbf{X})$ for $\eig(\mathbf{X}|\mathbf{Y})$.

\begin{theorem}
\label{thm:eigWmOrthogonal}
For any extension $\mathbf{X} \to \mathbf{Y}$ we have $\lp^2(\mathbf{X}) = \eig(\mathbf{X}|\mathbf{Y}) \oplus \wm(\mathbf{X}|\mathbf{Y})$.
\begin{proof}
Using the fact that the orbit of an eigenfunction is contained in a finite-rank $\lp^\infty(\mathbf{Y})$-module, one can show that every eigenfunction over $\mathbf{Y}$ is almost-periodic over $\mathbf{Y}$.
Thus $\eig(\mathbf{X}|\mathbf{Y}) \subset \wm(\mathbf{X}|\mathbf{Y})^\perp$ by Theorem~\ref{thm:ApWmOrthogonal}.

It remains to prove that $\eig(\mathbf{X}|\mathbf{Y})^\perp \subset \wm(\mathbf{X}|\mathbf{Y})$.
Fix $f$ in $\lp^2(\mathbf{X})$ orthogonal to $\eig(\mathbf{X}|\mathbf{Y})$.
For any $\phi$ in $\lp^\infty(\mathbf{X})$ we have
\begin{equation*}
\lim_{N \to \infty} \frac{1}{\haar(\Phi_N)} \int\limits_{\Phi_N} \! \int | \condex{\phi \cdot T^g f}{\mathbf{Y}}|^2 \intd\mu \intd\haar(g) = \langle f \conv H, f \rangle
\end{equation*}
as in the proof of Theorem~\ref{thm:ApWmOrthogonal}, where $H = \condex{\phi \otimes \phi}{\mathscr{I}}$.
Thus it suffices to prove that $f \conv H$ is in $\eig(\mathbf{X}|\mathbf{Y})$.
Let $\Psi_n$ and $\lambda_n$ be as in Theorem~\ref{joiningEigensections}.
Since $\{ \Psi_{n,y} \,:\, n \in \mathbb{N} \}$ spans the image of $H_y$ for almost-every $y$, it suffices to prove that each $\Psi_n$ is in $\eig(\mathbf{X}|\mathbf{Y})$.
To this end, fix $n$ in $\mathbb{N}$ and denote by $\theta(y)$ the multiplicity of the eigenvalue $\lambda_n(y)$ if $\lambda_n(y)$ is defined, and put $\theta(y) = 0$ otherwise.
Each of the functions $\lambda_m$ is measurable, so $\theta$ is too.
For each $k \in \mathbb{N}$ let $\Omega_k = \theta^{-1} (k)$.

We will show that the orbit of $1_{\Omega_k}(y) \Psi_n$ is a finite-rank $\lp^\infty(\mathbf{Y})$ module.
Fix $\gamma$ in $\Gamma$.
We have $T^\gamma H_{T^\gamma y} = H_y T^\gamma$ for almost-every $y$ because $H$ is $T \times T$ invariant.
Since $T^\gamma$ is unitary on almost every fiber, the operators $H_y$ and $H_{T^\gamma y}$ have the same spectrum, so each of the functions $\lambda_m$ is $T^\gamma$ invariant.
This implies $\Omega_k$ is $T^\gamma$-invariant.
Also
\begin{equation*}
(H \conv T^\gamma \Phi_n)_y = (T^\gamma (H \conv \Phi_n))_y = T^\gamma( H_{T^\gamma y} \conv \Phi_{n,T^\gamma y} ) = \lambda_n(y) (T^\gamma \Psi_n)_y
\end{equation*}
so in almost every fiber, the dimension the $\Gamma$ orbit of the square-integrable section $y \mapsto 1_{\Omega_k}(y) \Psi_{n,y}$ of $Y \total \mathfrak{H}$ is bounded by $k$.
Thus $y \mapsto 1_{\Omega_k}(y) \Psi_{n,y}$ corresponds to an eigenfunction.
Summing over $k$ proves that $\Psi_n$ is in $\eig(\mathbf{X}|\mathbf{Y})$ as desired.
\end{proof}
\end{theorem}

Combining Theorems~\ref{thm:ApWmOrthogonal} and \ref{thm:eigWmOrthogonal} yields the following result, a basic version of which will be used later.

\begin{corollary}
\label{cor:ApEigEqual}
For any extension $\mathbf{X} \to \mathbf{Y}$ we have $\ap(\mathbf{X}|\mathbf{Y}) = \eig(\mathbf{X}|\mathbf{Y})$.
\end{corollary}

\section{The Characteristic Factors}
\label{sec:characteristicFactors}

In this section we define the characteristic factors $\mathscr{C}_{k,i}$ associated to commuting, measurable actions $T_1,\dots,T_k$ and prove that
\begin{equation}
\label{eqn:characteristicFactorsStrongLimit}
\lim_{N \to \infty} \frac{1}{\haar(\Phi_N)} \int\limits_{\Phi_N} \prod_{i=1}^k T_k^g \cdots T_i^g f_i  - \prod_{i=1}^k T_k^g \cdots T_i^g \condex{f_i}{\mathscr{C}_{k,i}} \intd\haar(g) = 0
\end{equation}
in $\lp^2(X,\mathscr{B},\mu)$ for any $f_1,\dots,f_k$ in $\lp^\infty(X,\mathscr{B},\mu)$.

\begin{figure}
\begin{tikzpicture}
\matrix(a)[matrix of math nodes, row sep=2em, column sep=4em, text height=1.5ex, text depth=0.25ex]
{\mathscr{C}_{4,1} & \mathscr{C}_{4,2} & \mathscr{C}_{4,3} & \mathscr{C}_{4,4}\\ \mathscr{C}_{3,1} & \mathscr{C}_{3,2} & \mathscr{C}_{3,3} & \mathscr{C}_{3}\\ \mathscr{C}_{2,1} & \mathscr{C}_{2,2} & \mathscr{C}_2 & \\ \mathscr{C}_{1,1} & \mathscr{C}_1 & &\\};
\path[-] (a-1-1) edge node[left]{$T_4T_3T_2T_1$} (a-2-1);
\path[-] (a-1-2) edge node[left]{$T_4T_3T_2$} (a-2-2);
\path[-] (a-1-3) edge node[left]{$T_4T_3$} (a-2-3);
\path[-] (a-1-4) edge node[left]{$T_4$} (a-2-4);
\path[-] (a-2-1) edge node[left]{$T_3T_2T_1$} (a-3-1);
\path[-] (a-2-2) edge node[left]{$T_3T_2$} (a-3-2);
\path[-] (a-2-3) edge node[left]{$T_3$} (a-3-3);
\path[-] (a-3-1) edge node[left]{$T_2T_1$} (a-4-1);
\path[-] (a-3-2) edge node[left]{$T_2$} (a-4-2);
\end{tikzpicture}
\caption{The sub-$\sigma$-algebras $\mathscr{C}_{k,i}$ for $k \le 4$. A line indicates that the upper $\sigma$-algebra corresponds to the functions almost-periodic for the labeled action over the lower $\sigma$-algebra.}
\label{fig:CharacteristicFactorsSchematic}
\end{figure}
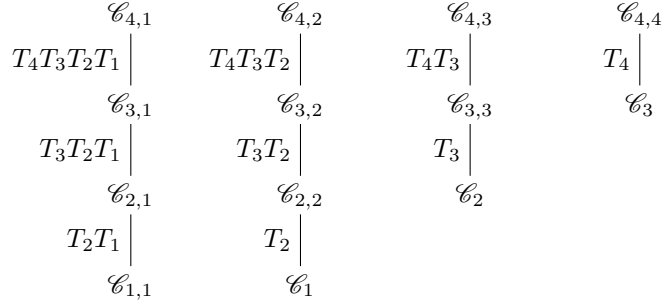

To define the $\mathscr{C}_{k,i}$ fix $k$ in $\mathbb{N}$ and let $T_1,\dots,T_k$ be commuting, measurable actions of $G$ on a probability space $(X,\mathscr{B},\mu)$.
Let $\mathscr{C}_{1,1}$ be the sub-$\sigma$-algebra of $T_1$ invariant sets.
It is invariant under all of the actions $T_2,\dots,T_k$ because they each commute with $T_1$.
Suppose by induction that for some $1 \le l \le k-1$ we have defined sub-$\sigma$-algebras $\mathscr{C}_{l,1},\dots,\mathscr{C}_{l,l}$ such that
\begin{enumerate}
\item for each $1 \le j \le l$, $\mathscr{C}_{l,j}$ is $T_l \cdots T_j$ invariant;
\item for each $1 \le j \le l$ and every $l+1 \le i \le k$, $\mathscr{C}_{l,j}$ is $T_i$ invariant.
\end{enumerate}
For each $1 \le j \le l$, let $\mathbf{Y}_j$ be the factor of $\mathbf{X}_j = (X,\mathscr{B},\mu,T_{l+1} \cdots T_j)$ corresponding to $\mathscr{C}_{l,j}$ and let $\mathscr{C}_{l+1,j}$ be the sub-$\sigma$-algebra of $\mathscr{B}$ corresponding to $\ap(\mathbf{X}_j|\mathbf{Y}_j)$.
It is invariant under $T_{l+1} \cdots T_j$ because it consists of $T_{l+1} \cdots T_j$ almost periodic functions, and (if $l < k-1$) it is $T_i$ invariant for all $l+2 \le i \le k$ because the actions commute.
Let $\mathbf{Y}_{l+1}$ be the factor of $\mathbf{X}_{l+1} = (X,\mathscr{B},\mu,T_{l+1})$ corresponding to $\mathscr{C}_{l,1} \vee \cdots \vee \mathscr{C}_{l,l}$ and let $\mathscr{C}_{l+1,l+1}$ be the sub-$\sigma$-algebra corresponding to $\ap(\mathbf{X}_{l+1}|\mathbf{Y}_{l+1})$.
It is $T_{l+1}$ invariant because it consists of the $T_{l+1}$ almost-periodic functions over $\mathbf{Y}_{l+1}$, and (if $l < k-1$) it is $T_i$ invariant for all $l+2 \le i \le k$ because the actions commute.
This concludes the inductive construction.
Figure \ref{fig:CharacteristicFactorsSchematic} shows how the $\mathscr{C}_{k,i}$ are related for $k \le 4$.
The remainder of this section constitutes a proof of the following theorem.

\begin{theorem}
\label{thm:characteristicFactors}
Let $T_1,\dots,T_k$ be commuting, measurable actions of $G$ on a separated, countably generated probability space $(X,\mathscr{B},\mu)$. Then
\begin{equation*}
\lim_{N \to \infty} \frac{1}{\haar(\Phi_N)} \int\limits_{\Phi_N} \prod_{i=1}^k T_k^g \cdots T_i^g f_i -  \prod_{i=1}^k T_k^g \cdots T_i^g \condex{f_i}{\mathscr{C}_{k,i}} \intd\haar(g) = 0
\end{equation*}
for any $f_i$ in $\lp^\infty(X,\mathscr{B},\mu)$.
\end{theorem}

Since the limit
\begin{equation*}
\lim_{N \to \infty} \frac{1}{\haar(\Phi_N)} \int\limits_{\Phi_N} \prod_{i=1}^k T_k^g \cdots T_i^g f_i \intd\haar(g)
\end{equation*}
is known to exist (see \cite{arXiv:1111.7292}) it suffices to prove that
\begin{equation}
\label{eqn:characteristicFactorsWeakLimit}
\begin{aligned}
 & \lim_{N \to \infty} \frac{1}{\haar(\Phi_N)} \int\limits_{\Phi_N} \! \int f_{k+1} \cdot \prod_{i=1}^k T_k^g \cdots T_i^g f_i \intd\mu \intd\haar(g)\\
= & \lim_{N \to \infty} \frac{1}{\haar(\Phi_N)} \int\limits_{\Phi_N} \! \int f_{k+1} \cdot \prod_{i=1}^k T_k^g \cdots T_i^g \condex{f_i}{\mathscr{C}_{k,i}} \intd\mu \intd\haar(g)
\end{aligned}
\end{equation}
for any $f_1,\dots,f_{k+1}$ in $\lp^\infty(X,\mathscr{B},\mu)$.
We will prove \eqref{eqn:characteristicFactorsWeakLimit} by induction on $k$.
The case $k = 1$ follows from the mean ergodic theorem: we have
\begin{equation*}
\lim_{N \to \infty} \frac{1}{\haar(\Phi_N)} \int\limits_{\Phi_N} \! \int T_1^g f_1 \cdot f_2 \intd\mu \intd\haar(g) = \int \condex{f_1}{\mathscr{C}_{1,1}} \cdot f_2 \intd\mu
\end{equation*}
by Proposition~\ref{meanErgodic}, which can be re-written as
\begin{equation*}
\int f_1 \otimes f_2 \intd\nu_1 = \int \condex{f_1}{\mathscr{C}_{1,1}} \otimes f_2 \intd\nu_1
\end{equation*}
where $\nu_1$ is the Furstenberg joining for the action $T_1$.
For the inductive step we need the following application of the van der Corput trick, which is a version of Lemma~4.7 in \cite{MR2599882}.

\begin{theorem}
\label{thm:liftingCharacteristicFactors}
Let $T_1,\dots,T_k$ be commuting, measurable actions of $G$ on a separated, countably generated probability space $(X,\mathscr{B},\mu)$.
Let $\nu_{k-1}$ be the Furstenberg joining of the actions $T_1,\dots,T_{k-1}$ and let $\nu_k$ be the Furstenberg joining of the actions $T_1,\dots,T_k$.
Suppose we have sub-$\sigma$-algebras $\mathscr{E}_1,\dots,\mathscr{E}_{k-1}$ with each $\mathscr{E}_i$ invariant under $T_{k-1} \cdots T_i$ and $T_k$ such that
\begin{equation*}
\int f_1 \otimes \cdots \otimes f_k \intd\nu_{k-1} = \int \condex{f_1}{\mathscr{E}_1} \otimes \cdots \otimes \condex{f_{k-1}}{\mathscr{E}_{k-1}} \otimes f_k \intd\nu_{k-1}
\end{equation*}
for all $f_1,\dots,f_k$ in $\lp^\infty(X,\mathscr{B},\mu)$.
Put $\mathscr{E}_k = \mathscr{E}_1 \vee \cdots \vee \mathscr{E}_{k-1}$.
Then
\begin{equation}
\label{eqn:characteristicInductionGoal}
\int f_1 \otimes \cdots \otimes f_{k+1} \intd\nu_k = \int \condex{f_1}{\mathscr{F}_1} \otimes \cdots \otimes \condex{f_k}{\mathscr{F}_k} \otimes f_{k+1} \intd\nu_k
\end{equation}
for all $f_1,\dots,f_{k+1}$ in $\lp^\infty(X,\mathscr{B},\mu)$ where, for each $1 \le i \le k$ the sub-$\sigma$-algebra $\mathscr{F}_i$ corresponds to the functions that are $T_k \cdots T_i$ almost-periodic over $\mathscr{E}_i$.
\begin{proof}
Fix $f_1,\dots,f_{k+1}$ in $\lp^\infty(X,\mathscr{B},\mu)$ with $\nbar f_i \nbar_\infty \le 1$ for all $1 \le i \le k+1$.
Since $\mathscr{E}_i$ is $T_{k-1} \cdots T_i$ invariant and contained in $\mathscr{E}_k$ for each $1 \le i \le k-1$, we have
\begin{equation*}
T_{k-1}^g \cdots T_i^g \condex{f_i}{\mathscr{E}_i} = \condex{T_{k-1}^g \cdots T_i^g \condex{f_i}{\mathscr{E}_i}}{\mathscr{E}_k}
\end{equation*}
for each $1 \le i \le k-1$.
Thus we can re-write our assumption as
\begin{equation}
\label{eqn:characteristicInductionStep}
\begin{aligned}
\int f_1 \otimes \cdots \otimes f_k \intd\nu_{k-1} = \int \condex{f_1}{\mathscr{E}_1} \otimes \cdots \otimes \condex{f_k}{\mathscr{E}_k} \intd\nu_{k-1}
\end{aligned}
\end{equation}
using \eqref{FurstenbergJoining}.
We proceed by applying the van der Corput trick to the sequence
\begin{equation*}
u(g) = \prod_{i=1}^k T_k^g \cdots T_i^g f_i
\end{equation*}
in $\lp^2(X,\mathscr{B},\mu)$.
From \eqref{FurstenbergJoining} and \eqref{eqn:characteristicInductionStep} we see that
\begin{align*}
\lim_{N \to \infty} \frac{1}{\haar(\Phi_N)} \int\limits_{\Phi_N} \! \langle u(hg), & u(lg) \rangle \intd\haar(g) = \int \bigotimes_{i=1}^k \left( T_k^h \cdots T_i^h f_i \cdot T_k^l \cdots T_i^l f_i \right) \intd\nu_{k-1}\\
 = & \int \bigotimes_{i=1}^k \condex{ T_k^h \cdots T_i^h f_i \cdot T_k^l \cdots T_i^l f_i}{\mathscr{E}_i} \intd\nu_{k-1}
\end{align*}
for any $h,l \in G$.
Using \eqref{FurstenbergJoining} once more yields
\begin{equation*}
\lim_{N \to \infty} \frac{1}{\haar(\Phi_N)} \int\limits_{\Phi_N} \! \langle u(hg), u(lg) \rangle \intd\haar(g) \le \nbar \condex{T_k^h \cdots T_i^h f_i \cdot T_k^l \cdots T_i^l f_i}{\mathscr{E}_i} \nbar
\end{equation*}
for each $1 \le i \le k$, the norm taken in $\lp^2(X,\mathscr{B},\mu)$.
Let $\mathbf{X}_i$ be the system $(X,\mathscr{B},\mu,T_k \cdots T_i)$ and let $\mathbf{Y}_i$ be a factor corresponding to $\mathscr{E}_i$.
Let $\mathscr{I}_i$ be the sub-$\sigma$-algebra of $T_k \cdots T_i \times T_k \cdots T_i$ invariant sets in the relatively independent joining $\mathbf{X}_i \times_{\mathbf{Y}_i} \mathbf{X}_i$.
We have
\begin{align*}
 & \limsup_{H \to \infty} \frac{1}{\haar(\Phi_H)^2} \int\limits_{\Phi_H} \! \int\limits_{\Phi_H} \nbar \condex{T_k^h \cdots T_i^h f_i \cdot T_k^l \cdots T_i^l f_i}{\mathscr{E}_i} \nbar \intd\haar(h) \intd\haar(l)\\
\le & \lim_{H \to \infty} \bnbar \frac{1}{\haar(\Phi_H)} \int\limits_{\Phi_H} (T_k \cdots T_i \times T_k \cdots T_i)^h (f_i \otimes f_i) \intd\haar(h) \bnbar = \nbar \condex{f_i \otimes f_i}{\mathscr{I}_i} \nbar
\end{align*}
in $\lp^2(\mathbf{X}_i \times_{\mathbf{Y}_i} \mathbf{X}_i)$ by Cauchy-Schwarz and the mean ergodic theorem.
By Corollary~\ref{invariantWmIsZero} the conditional expectation $\condex{f_i \otimes f_i}{\mathscr{I}_i}$ will be zero if $f_i$ is orthogonal to $\ap(\mathbf{X}_i|\mathbf{Y}_i)$.
Since $1 \le i \le k$ was arbitrary, \eqref{eqn:characteristicInductionGoal} follows from the van der Corput trick.
\end{proof}
\end{theorem}

Taking $\mathscr{E}_i = \mathscr{C}_{k-1,i}$ in the preceding theorem proves \eqref{eqn:characteristicFactorsWeakLimit} and concludes the proof of Theorem~\ref{thm:characteristicFactors}.
We conclude this section with another application of the van der Corput trick that is sometimes useful.

\begin{theorem}
\label{UsingVanDerCorput}
Let $T_1,\dots,T_k$ be commuting, measurable actions of $G$ on a separated, countably generated probability space $(X,\mathscr{B},\mu)$.
Let $\nu_{k-1}$ be the Furstenberg joining of the actions $T_1,\dots,T_{k-1}$ and let $\nu_k$ be the Furstenberg joining of the actions $T_1,\dots,T_k$.
Let $\mathscr{I}_k$ denote the sub-$\sigma$-algebra of $\mathscr{B}^k$ consisting of
\begin{equation*}
T_k T_{k-1} \cdots T_1 \times \cdots \times T_k T_{k-1} \times T_k
\end{equation*}
invariant sets.
If $\condex{f_1 \otimes \cdots \otimes f_k}{\mathscr{I}_k} = 0$ in $\lp^2(X^k,\mathscr{B}^k,\nu_{k-1})$ for some $f_1,\dots,f_k$ in $\lp^\infty(X,\mathscr{B},\mu)$ then
\begin{equation*}
\int f_1 \otimes \cdots \otimes f_k \otimes f_{k+1} \intd\nu_k = 0
\end{equation*}
for all $f_{k+1}$ in $\lp^2(X,\mathscr{B},\mu)$.
\begin{proof}
Fix $f_1,\dots,f_k$ in $\lp^\infty(X,\mathscr{B},\mu)$ satisfying $\condex{f_1 \otimes \cdots \otimes f_k}{\mathscr{I}_k} = 0$.
Applying the van der Corput trick as in Theorem~\ref{thm:liftingCharacteristicFactors} gives
\begin{equation*}
\lim_{N \to \infty} \frac{1}{\haar(\Phi_N)} \int\limits_{\Phi_N} \! \langle u(hg), u(lg) \rangle \intd\haar(g) = \int \bigotimes_{i=1}^k \left( T_k^h \cdots T_i^h f_i \cdot T_k^l \cdots T_i^l f_i \right) \intd\nu_{k-1}
\end{equation*}
for any $h,l \in G$.
From this we get
\begin{align*}
 & \lim_{H \to \infty} \frac{1}{\haar(\Phi_H)^2} \int\limits_{\Phi_H} \!\int\limits_{\Phi_H} \lim_{N \to \infty} \frac{1}{\haar(\Phi_N)} \int\limits_{\Phi_N} \! \langle u(hg), u(lg) \rangle \intd\haar(g) \intd\haar(h) \intd\haar(l)\displaybreak[2]\\
= & \lim_{H \to \infty} \bnbar \frac{1}{\haar(\Phi_H)} \int\limits_{\Phi_H} \bigotimes_{i=1}^k \, T_k^h \cdots T_i^h f_i \intd\haar(h) \bnbar^2 = \int \condex{f_1 \otimes \cdots \otimes f_k}{\mathscr{I}_k}^2 \intd\nu_{k-1}
\end{align*}
where the norm is determined by $\nu_{k-1}$ and the last equality follows from Proposition \ref{meanErgodic}.
The conclusion follows from the van der Corput trick and the fact that strong convergence implies weak convergence.
\end{proof}
\end{theorem}

\section{Lifting Positivity}
\label{sec:liftingPositivity}

In this section we prove a technical result, based on Theorem~9.1 in \cite{MR670131}, that allows us to lift multiple recurrence from one level of Figure~\ref{fig:CharacteristicFactorsSchematic} to the next, provided the sub-$\sigma$-algebras in the lower level are all equal.
In the next section we will use this to prove some multiple recurrence results.

\begin{theorem}
\label{thm:mainPositivityResult}
Let $T_1,\dots,T_k$ be commuting, measurable actions of $G$ on a separated, countably-generated probability space $(X,\mathscr{B},\mu)$.
Let $\mathscr{D}$ be a sub-$\sigma$-algebra that is $T_k \cdots T_i$ invariant for all $1 \le i \le k$.
Suppose that
\begin{equation*}
\liminf_{N \to \infty} \frac{1}{\haar(\Phi_N)} \int\limits_{\Phi_N} \! \int f \prod_{i=1}^k T_k^g \cdots T_i^g \condex{f}{\mathscr{D}} \intd\mu \intd\haar(g) > 0
\end{equation*}
for any $f > 0$ in $\lp^\infty(X,\mathscr{B},\mu)$.
For each $1 \le i \le k$, let $\mathscr{E}_i$ be a sub-$\sigma$-algebra of $\mathscr{B}$ that is $T_k \cdots T_i$ invariant, and suppose that $\mathscr{E}_i \to \mathscr{D}$ is $T_k \cdots T_i$ almost-periodic.
Then
\begin{equation*}
\liminf_{N \to \infty} \frac{1}{\haar(\Phi_N)} \int\limits_{\Phi_N} \! \int f \prod_{i=1}^k T_k^g \cdots T_i^g \condex{f}{\mathscr{E}_i} \intd\mu \intd\haar(g) > 0
\end{equation*}
for any $f > 0$ in $\lp^\infty(X,\mathscr{B},\mu)$.
\begin{proof}
It suffices to prove the theorem when $f = 1_B$ for some set $B \in \mathscr{B}$ having positive measure.
Let $\mu_x$ be the disintegration of $\mu$ over $\mathscr{D}$.
For each $1 \le i \le k$ write $f_i$ for $\condex{1_B}{\mathscr{E}_i}$.
From
\begin{equation*}
\mu(B \cap \{ f_i = 0 \}) = \int f_i \cdot 1_{\{ f_i = 0 \}} \intd\mu = 0
\end{equation*}
it follows that $f_i$ is positive on almost all of $B$.
Thus we can find a set $D_1$ in $\mathscr{D}$ with positive measure and some $\alpha > 0$ such that
\begin{equation}
\label{positiveProduct}
\int f \cdot f_1 \cdots f_k \intd\mu_x > \alpha
\end{equation}
for all $x$ in $D_1$.
Fix $\epsilon = \alpha/4k$.

For any $x \in X$ and any non-empty subset $F$ of $\Gamma$ define
\begin{equation*}
\mathcal{L}(x,F) = \{ (T_k^a \cdots T_1^a f_1, \dots, T_k^a f_k) \,:\, a \in F \} \subset \lp^2(X,\mathscr{B},\mu_x)^k
\end{equation*}
and equip it with the max norm coming from $\nbar \cdot \nbar_x$ on the constituents.

\begin{claim}
There is a subset $D_2$ of $D_1$ with positive measure such that $\mathcal{L}(x,\Gamma)$ is totally bounded for each $x \in D_2$.
\end{claim}
\begin{proof}
For each $j$ in $\mathbb{N}$ put $\epsilon_j = \mu(D_1)/2^{j+k+1}$.
Since each $f_i$ is $T_k \cdots T_i$ almost-periodic over $\mathscr{D}$ one can find finite subsets $\Xi^i_j$ of $\lp^\infty(X,\mathscr{B},\mu)$ and subsets $E^i_j$ with measure at least $1 - \epsilon_j$ such that for each $\gamma$ in $\Gamma$ we have
\begin{align*}
\min \{ \nbar T_k^\gamma \cdots T_i^\gamma f_i - \xi \nbar_x \,:\, \xi \in \Xi^i_j \} &< \epsilon_j
\end{align*}
for every $x$ in $E^i_j$.
Put
\begin{equation*}
D_2 = D_1 \Big\backslash \bigcup_{j=1}^\infty E^1_j \cup \cdots \cup E^k_j
\end{equation*}
and note that $\mu(D_2) \ge \mu(D_1)/2$.
\end{proof}

We will be interested in separated subsets of $\mathcal{L}(x,\Gamma)$ so define
\begin{equation*}
\sep(F,t) := \bigcap_{a \in F} \bigcap_{b \in F}^{a \ne b} \bigcup_{i=1}^k \{ x \in X \,:\, \nbar T_k^a \cdots T_i^a f_i - T_k^b \cdots T_i^b f_i \nbar_x > t \}
\end{equation*}
for any finite, non-empty subset $F$ of $\Gamma$ and any positive $t$.
It belongs to $\mathscr{D}$ and when $F$ is a singleton it is all of $X$.
The fact that $\mathcal{L}(x,\Gamma)$ is totally bounded whenever $x \in D_2$ implies that there is a bound on the cardinality of the finite sets $F$ for which $x$ belongs to $\sep(F,\epsilon)$.
Thus the $\mathscr{D}$ measurable sets
\begin{equation*}
Q(F) := \sep(F,\epsilon) \Big\backslash \bigcup \, \{ \sep(E,\epsilon) \,:\, E \subset \Gamma \text{ with } |F| < |E| < \infty \}
\end{equation*}
cover almost all of $D_2$ as $F$ runs through the finite subsets of $\Gamma$ and we can fix a finite, non-empty subset $F$ of $\Gamma$ such that $Q(F) \cap D_2$ has positive measure.
For each $x$ in $Q(F) \cap D_2$ we can find some $n \in \mathbb{N}$ with the property that $x \in \sep(F,\epsilon + 1/n)$ because of the strict inequalities and finite number of conditions in the definition of $\sep(F,\epsilon)$.
Thus we can find some $\eta > 0$ with the property that $Q(F) \cap \sep(F,\epsilon + \eta) \cap D_2$ has positive measure.
Define a function $\Psi$ by
\begin{align*}
\Psi : Q(F) \cap \sep(F,\epsilon + \eta) \cap D_2 &\to [0,2]^{F \times F \times \{1,\dots,k \} }\\
\Psi(x) : (a,b,i) &\mapsto \nbar T_k^a \cdots T_i^a f_i - T_k^b \cdots T_i^b f_i \nbar_x
\end{align*}
and partition $[0,2]^{F \times F \times \{1,\dots,k\}}$ into cubes of side length $\eta/2$.
Since $\Psi$ is measurable we can find a cell $D$ in the pull-back partition that has positive measure.
Now $D$ belongs to $\mathscr{D}$, so by hypothesis
\begin{equation*}
\liminf_{N \to \infty} \frac{1}{\haar(\Phi_N)} \int\limits_{\Phi_N} \! \int \prod_{i=1}^k T_k^g \cdots T_i^g 1_D \cdot 1_D \intd\mu \intd\haar(g) > 0
\end{equation*}
and thus there exists $\zeta > 0$ and a subset $\Delta$ of $G$ with positive lower density such that
\begin{equation*}
\int \prod_{i=1}^k T_k^g \cdots T_i^g 1_D \cdot 1_D\intd\mu > \zeta
\end{equation*}
for any $g$ in $\Delta$.

\begin{claim}
For any $g \in \Delta$ there is a subset $E_g$ of $(T_k^g \cdots T_1^g)^{-1} D \cap \cdots \cap (T_k^g)^{-1} D \cap D$ with measure at least $\zeta/2$ such that for any $x \in E_g$ one can find $b \in F$ satisfying $\nbar T_k^{bg} \cdots T_i^{bg} f_i - f_i \nbar_x < 2\epsilon$ for every $1 \le i \le k$.
\end{claim}
\begin{proof}
Fix $g \in \Delta$.
Since $\Gamma$ is dense and each $T_i^g$ is unitary we can find $\gamma$ in $\Gamma$ such that
\begin{equation}
\label{densityOfSubgroupThing}
\nbar T_k^{ag} \cdots T_i^{ag} f_i - T_k^{a\gamma} \cdots T_i^{a\gamma} f_i\nbar^2 \le \min \{ \eta^2 \zeta/2^{2k+2}|F|, \epsilon^2 \}
\end{equation}
for all $a \in F$ and all $1 \le i \le k$.
It follows from Chebyshev's inequality that there is a subset $E_g$ of $(T_k^g \cdots T_1^g)^{-1} D \cap \cdots \cap (T_k^g)^{-1} D \cap D$ with $\mu(E_g) \ge \zeta/2$ such that
\begin{equation}
\label{stupidName}
\nbar T_k^{ag} \cdots T_i^{ag} f_i - T_k^{a\gamma} \cdots T_i^{a\gamma} f_i\nbar_x \le \eta/4
\end{equation}
for all $a \in F$, all $1 \le i \le k$ and all $x \in E_g$.

If $1 \in F\gamma$ then the claim follows immediately from \eqref{densityOfSubgroupThing}, so assume otherwise.
In this case the subset $F\gamma \cup \{1\}$ of the subgroup $\Gamma$ has cardinality strictly larger than $F$ so $x$ does not belong to $\sep(F\gamma \cup \{ 1 \},\epsilon)$.
Thus we can find $\alpha \ne \beta$ in $F\gamma \cup \{ 1 \}$ such that
\begin{equation}
\label{ChoosingFinitelyWisely}
\nbar T_k^\alpha \cdots T_i^\alpha f_i - T_k^\beta \cdots T_i^\beta f_i \nbar_x \le \epsilon
\end{equation}
for all $1 \le i \le k$, and the proof will be concluded if we can show that one of $\alpha$ or $\beta$ must be 1.
Fix $a \ne b$ in $F$.
(If $|F| = 1$ then one of $\alpha$ or $\beta$ must be 1.)
That $x$ belongs to $\sep(F,\epsilon + \eta)$ tells us
\begin{equation*}
\nbar T_k^a \cdots T_i^a f_i - T_k^b \cdots T_i^b f_i \nbar_x > \epsilon + \eta
\end{equation*}
holds for some $1 \le i \le k$.
Since $T_k^g \cdots T_i^g x$ belongs to $D$ we must have
\begin{equation*}
\nbar T_k^a \cdots T_i^a f_i - T_k^b \cdots T_i^b f_i \nbar_{T_k^g \cdots T_i^g x} > \epsilon +\eta/2
\end{equation*}
because the function $x \mapsto \Psi(x)(a,b,i)$ takes values in an interval of length at most $\eta/2$.
Now $a \ne b$ in $F$ were arbitrary so, combined with \eqref{stupidName}, this forces one of $\alpha$ or $\beta$ to be 1 as otherwise \eqref{ChoosingFinitelyWisely} is contradicted.
\end{proof}

We can now finish the proof.
Fix $g \in \Delta$ and let $E_g$ be as in the claim.
For any $x \in E_g$ we can find some $b \in F$ such that $\nbar T_k^{bg} \cdots T_i^{bg} f_i - f_i \nbar_x \le 2\epsilon$ for all $1 \le i \le k$.
Thus
\begin{equation*}
\int f \cdot \prod_{i=1}^k T_k^{bg} \cdots T_i^{bg} f_i \intd\mu_x \ge \alpha - 2k\epsilon = \frac{\alpha}{2}
\end{equation*}
for any $x$ in the subset $E_g$ of $D$.
Summing over $b \in F$ on the left hand side weakens the inequality and removes the dependence of $b$ on $x$ and $g$.
This allows us to integrate over $E_g$, obtaining
\begin{equation*}
\sum_{b \in F} \int f \cdot \prod_{i=1}^k T_k^{bg} \cdots T_i^{bg} f_i \intd\mu \ge \frac{\zeta\alpha}{4}
\end{equation*}
which, after averaging over $\Delta$ using the F{\o}lner sequence $\Phi$, gives
\begin{equation*}
\liminf_{N \to \infty} \frac{1}{\haar(\Phi_N)} \int\limits_{\Phi_N} \! \int f \prod_{i=1}^k T_k^g \cdots T_i^g f_i \intd\mu \intd\haar(g) \ge \liminf_{N \to \infty} \frac{\haar(\Delta \cap \Phi_N)}{\haar(\Phi_N)} \cdot \frac{\zeta \alpha}{4|F|}
\end{equation*}
concluding the proof.
\end{proof}
\end{theorem}

\section{Recurrence Results}
\label{sec:recurrenceResults}

Bergelson's conjecture states that for any commuting actions $T_1,\dots,T_k$ of $G$ on a separated, countably generated probability space $(X,\mathscr{B},\mu)$ we have
\begin{equation}
\label{eqn:amenableSz}
\liminf_{N \to \infty} \frac{1}{\haar(\Phi_N)} \int\limits_{\Phi_N} \! \int 1_B \cdot \prod_{i=1}^k T_k \cdots T_i^g 1_B \intd\mu \intd \haar(g) > 0
\end{equation}
for every $B$ in $\mathscr{B}$ with positive measure.
In this section we verify this conjecture when $k = 2$ without additional assumptions, and when $k = 3$ assuming $T_1,T_2$ and $T_2T_1$ are ergodic.
The $k = 2$ case was previously obtained for countable, amenable groups in \cite{MR1481813}.

\begin{theorem}
\label{thm:amenableRoth}
Let $T_1,T_2$ be commuting, measurable actions of $G$ on a separated, countably generated probability space $(X,\mathscr{B},\mu)$.
Then
\begin{equation}
\label{eqn:amenableRoth}
\liminf_{N \to \infty} \frac{1}{\haar(\Phi_N)} \int\limits_{\Phi_N} \! \int f \cdot T_2^g T_1^g f \cdot T_2^g f \intd\mu \intd\haar(g) > 0
\end{equation}
for any $f > 0$ in $\lp^\infty(X,\mathscr{B},\mu)$.
\begin{proof}
By Theorem \ref{thm:characteristicFactors} it suffices to prove that
\begin{equation}
\label{eqn:RothCharacteristicPositivity}
\liminf_{N \to \infty} \frac{1}{\haar(\Phi_N)} \int\limits_{\Phi_N} \! \int f \cdot T_2^g T_1^g \condex{f}{\mathscr{C}_{2,1}} \cdot T_2^g \condex{f}{\mathscr{C}_{2,2}} \intd \mu \intd \haar(g) > 0
\end{equation}
for all $f > 0$ in $\lp^\infty(X,\mathscr{B},\mu)$.
If $f$ is of the form $1_B$ for some $B \in \mathscr{C}_1$ with $\mu(B) > 0$ then $T_2 T_1 f \cdot T_2 f = T_2 f$ so in this case \eqref{eqn:amenableRoth} follows from the mean ergodic theorem.
Thus we have \eqref{eqn:amenableRoth} whenever $f$ is $\mathscr{C}_1$ measurable.
Applying Theorem~\ref{thm:mainPositivityResult} with $\mathscr{D} = \mathscr{C}_1$, $\mathscr{E}_1 = \mathscr{C}_{2,1}$ and $\mathscr{E}_2 = \mathscr{C}_{2,2}$ yields \eqref{eqn:RothCharacteristicPositivity}.
\end{proof}
\end{theorem}

When $k = 3$ we cannot use Theorem~\ref{thm:mainPositivityResult} to prove \eqref{eqn:amenableSz} because the sub-$\sigma$-algebras $\mathscr{C}_{2,1}$ and $\mathscr{C}_{2,2}$ need not agree and because the behavior of $T_3$ with respect to the extensions $\mathscr{C}_{2,i} \to \mathscr{C}_{1,1}$ is unknown.
However, if $T_1$ is ergodic then $\mathscr{C}_{1,1}$ is trivial and the sub-$\sigma$-algebras $\mathscr{C}_{2,1}$ and $\mathscr{C}_{2,2}$ consist of functions that are almost-periodic for $T_2T_1$ and $T_2$ respectively over the trivial factor.
We will prove below that if $T_2$ and $T_2 T_1$ are ergodic any function almost-periodic for $T_2$ or $T_2T_1$ over the trivial factor is necessarily almost periodic for $T_3$ over the trivial factor.
This leads to a description of characteristic factors that allow us, under the aforementioned ergodicity assumptions, to prove Bergelson's conjecture when $k = 3$.

Given a system $\mathbf{X}$, recall that $f$ in $\lp^2(\mathbf{X})$ is an \define{eigenfunction} of $\mathbf{X}$ if its $T$-orbit is contained in a $T$-invariant, finite-dimensional subspace of $\lp^2(\mathbf{X})$.
In other words $f$ is an eigenfunction of $T$ if its orbit is contained in a finite-dimensional sub-representation of $\lp^2(\mathbf{X})$.
Denote by $\eig(\mathbf{X})$ or $\eig(T)$ the closure of the subspace of $\lp^2(\mathbf{X})$ spanned by the eigenfunctions of $\mathbf{X}$.

\begin{proposition}
\label{prop:ergodicEigenfunctionsForCommuting}
Let $S_1$ and $S_2$ be commuting actions of $G$ on a separated, countably generated probability space $(X,\mathscr{B},\mu)$.
If $S_2$ is ergodic then $\eig(S_2) \subset \eig(S_1)$.
\begin{proof}
Let $f$ be an eigenfunction of $S_2$ and let $\mathscr{M}$ be an $S_2$-invariant, finite-dimensional subspace of $\lp^2(X,\mathscr{B},\mu)$ containing the orbit of $f$.
Without loss of generality, we can assume $\mathscr{M}$ is irreducible.
Let $\mathscr{N}$ be the closed subspace of $\lp^2(X,\mathscr{B},\mu)$ spanned by the sub-representations of $G$ on $\lp^2(X,\mathscr{B},\mu)$ induced by $S_2$ that are equivalent to $\mathscr{M}$.
By Proposition~1.4 in \cite{MR961735} the multiplicity of $\mathscr{M}$ in $\lp^2(X,\mathscr{B},\mu)$ is bounded by its dimension, so $\mathscr{N}$ is finite-dimensional.
Fix $g \in G$ and put $\mathscr{M}_g = \{ S_1^g f \,:\, f \in \mathscr{M} \}$.
Since $S_1$ and $S_2$ commute the representations of $G$ on $\mathscr{M}$ and $\mathscr{M}_g$ determined by $S_2$ are equivalent.
Thus $\mathscr{M}_g \subset \mathscr{N}$.
This implies $S_1^g f \in \mathscr{N}$ for all $g \in G$, so $f$ is contained in $\eig(S_1)$.
\end{proof}
\end{proposition}

\begin{proposition}
\label{prop:ergodicEigenfunctionsForProduct}
Let $S_1$ and $S_2$ be commuting actions of $G$ on a separated, countably generated probability space $(X,\mathscr{B},\mu)$. If $S_2$ is ergodic then $\eig(S_2) \subset \eig(S_2 S_1)$.
\begin{proof}
Let $f$ be an eigenfunction of $S_2$.
Form $\mathscr{M}$ and $\mathscr{N}$ as in the proof of Proposition~\ref{prop:ergodicEigenfunctionsForCommuting}.
Fix $g \in G$.
Put $\mathscr{M}_g = S_2^g S_1^g \mathscr{M}$.
We have $\mathscr{M}_g = S_1^g \mathscr{M}$, which is equivalent to $\mathscr{M}$ and therefore contained in $\mathscr{N}$, as desired.
\end{proof}
\end{proposition}

We can now give a proof of Bergelson's conjecture when $k = 3$ and the actions $T_1,T_2$ and $T_2T_1$ are all ergodic.

\begin{theorem}
\label{thm:amenableFourTermSz}
Let $T_1,T_2,T_3$ be commuting, measurable actions of $G$ on a separated, countably generated probability space $(X,\mathscr{B},\mu)$.
Suppose that the actions $T_1, T_2$ and $T_2T_1$ are ergodic.
Then
\begin{equation}
\label{eqn:fourTermSz}
\liminf_{N \to \infty} \frac{1}{\haar(\Phi_N)} \int\limits_{\Phi_N} \! \int f \cdot T_3^g T_2^g T_1^g f \cdot T_3^g T_2^g f \cdot T_3^g f \intd\mu \intd\haar(g) > 0
\end{equation}
for any $f > 0$ in $\lp^\infty(X,\mathscr{B},\mu)$.
\begin{proof}
Ergodicity of $T_1$ means $\mathscr{C}_{1,1}$ is trivial, so $\mathscr{C}_{2,1}$ and $\mathscr{C}_{2,2}$ correspond to the functions that are $T_2T_1$ and $T_2$ almost-periodic over the trivial factor respectively.
Let $\mathscr{D}$ be the sub-$\sigma$-algebra of $\mathscr{B}$ corresponding to the functions that are almost-periodic for $T_3$ over the trivial factor.
Combining Corollary~\ref{cor:ApEigEqual} with Propositions~\ref{prop:ergodicEigenfunctionsForCommuting} and \ref{prop:ergodicEigenfunctionsForProduct} gives $\mathscr{C}_{2,2} \subset \mathscr{C}_{2,1} \subset \mathscr{D}$.
This implies any $\mathscr{C}_{2,1}$ measurable function $f$ is almost-periodic for both $T_2T_1$ and $T_3$, so $f \in \eig(T_3T_2T_1)$.

We begin by showing that
\begin{equation*}
\liminf_{N \to \infty} \frac{1}{\haar(\Phi_N)} \int\limits_{\Phi_N} \! \int T_3^g T_2^g T_1^g f \cdot T_3^g T_2^g f \cdot T_3^g f \cdot f \intd\haar(g) \intd\mu > 0
\end{equation*}
whenever $f = 1_B$ is $\mathscr{C}_{2,1}$ measurable and $\mu(B) > 0$.
Put $\epsilon = \mu(B)^2/6$.
Since $\mathscr{C}_{2,1} \subset \mathscr{D}$ the set $\Omega_1$ of $g$ in $G$ for which $\nbar T_3^g 1_B - 1_B \nbar \le \epsilon$ and $\nbar T_3^g T_2^g T_1^g 1_B - 1_B \nbar \le \epsilon$ is a measurable $\ip^*$ subset of $G$.
By the argument on page~50 of \cite{MR1411215} the set $\Omega_2$ consisting of those $g$ for which $\mu(B \cap (T_3^g T_2^g)^{-1} B) \ge \mu(B)^2 - \epsilon$ is $\ip^*$.
It is also measurable, so the intersection $\Omega = \Omega_1 \cap \Omega_2$ is measurable and $\ip^*$.
We have
\begin{equation*}
\begin{aligned}
 & \liminf_{N \to \infty} \frac{1}{\haar(\Phi_N)} \int\limits_{\Phi_N} \int T_3^g T_2^g T_1^g 1_B \cdot T_3^g T_2^g 1_B \cdot T_3^g 1_B \cdot 1_B \intd\mu \intd\haar(g)\\
\ge & \liminf_{N \to \infty} \frac{1}{\haar(\Phi_N)} \int\limits_{\Phi_N} 1_\Omega(g) \left( \int T_3^g T_2^g 1_B \cdot 1_B \intd\mu - 2\epsilon \right) \intd\haar(g)\\
\ge & \liminf_{N \to \infty} \frac{1}{\haar(\Phi_N)} \int\limits_{\Phi_N} 1_\Omega(g) \intd\haar(g) \cdot \left( \mu(B)^2 - 3\epsilon \right) = \frac{\lowerdens(\Omega) \mu(B)^2}{2} > 0
\end{aligned}
\end{equation*}
because every $\ip^*$ set has positive lower density.

The fact that $\mathscr{C}_{2,2} \subset \mathscr{C}_{2,1}$ implies
\begin{equation*}
\int f_1 \otimes f_2 \otimes f_3 \intd\nu_2 = \int \condex{f_1}{\mathscr{C}_{2,1}} \otimes \condex{f_2}{\mathscr{C}_{2,1}} \otimes \condex{f_3}{\mathscr{C}_{2,1}} \intd\nu_2
\end{equation*}
for all $f_1,f_2,f_3$ in $\lp^\infty(X,\mathscr{B},\mu)$ where $\nu_2$ is the Furstenberg joining for the actions $T_1$ and $T_2$.
Applying Theorem~\ref{thm:liftingCharacteristicFactors} with $\mathscr{D}_{2,i} = \mathscr{C}_{2,1}$ gives sub-$\sigma$-algebras $\mathscr{E}_{3,i}$ that are characteristic for \eqref{eqn:fourTermSz}.
Moreover $\mathscr{E}_{3,i} \to \mathscr{C}_{2,1}$ is compact for $T_3\cdots T_i$. 
Finally, using Theorem~\ref{thm:mainPositivityResult} with $k = 3$, $\mathscr{D} = \mathscr{C}_{2,1}$ and $\mathscr{E}_i = \mathscr{E}_{3,i}$ yields \eqref{eqn:fourTermSz}.
\end{proof}
\end{theorem}

Given commuting, measurable actions $T_1,\dots,T_k$ of $G$ on a separated, countably generated probability space $(X,\mathscr{B},\mu)$ define
\begin{equation*}
\ret_k(B) = \{ g \in G \,:\, \mu(B \cap (T_k^g \cdots T_1^g)^{-1} B \cap \cdots \cap (T_k^g)^{-1} B) > 0 \}
\end{equation*}
for any $B$ in $\mathscr{B}$.
We say that a subset $R$ of $G$ is \define{syndetic} if there is a compact set $F \subset G$ such that $FS = G$.
The following result, based on the argument on page 1199 in \cite{MR1481813}, shows that $\ret_k(B)$ is syndetic whenever it has positive lower density with respect to any F\o{}lner sequence.

\begin{lemma}
Let $R$ be a measurable subset of $G$ that has positive lower density with respect to every F\o{}lner sequence. Then $R$ is syndetic.
\begin{proof}
Suppose $R$ is not syndetic.
Then for every compact subset $F$ of $G$ we have $FR \ne G$.
For each $N \in \mathbb{N}$ choose $h_N$ from $G \backslash \Phi_N^{-1} R$.
Then $\Phi_N h_N \cap R$ is empty for all $N$.
However, $N \mapsto \Phi_N h_N$ is a left F\o{}lner sequence (because the modular function is everywhere positive) so $R$ must have positive lower density with respect to it, giving the desired contradiction.
\end{proof}
\end{lemma}

It now follows immediately from Theorems~\ref{thm:amenableRoth} and \ref{thm:amenableFourTermSz} that $R_2$ is always syndetic and that $R_3$ is syndetic if $T_1,T_2$ and $T_2T_1$ are ergodic.

\section{Further Results}
\label{sec:furtherResults}

Recall that a system $\mathbf{X}$ is said to be \define{weakly mixing} if the only functions almost-periodic over the trivial factor consisting of one point are the constant functions.
An immediate consequence of Theorem \ref{thm:characteristicFactors} is that the limit as $N \to \infty$ of \eqref{eqn:amenableCommutingAverage} is constant whenever all of the systems $(X,\mathscr{B},\mu,T_j \cdots T_i)$ are weakly mixing.
This is because all the $\mathscr{C}_{j,i}$ are trivial in that case.
In fact, we have a short proof of the following result.

\begin{theorem}[2.4 in \cite{MR961735}]
\label{thm:jointErgodicityImpliesConstantLimit}
Let $T_1,\dots,T_k$ be commuting, measurable actions of $G$ on a separated, countably generated probability space $(X,\mathscr{B},\mu)$.
If each of the actions
\begin{equation*}
T_1,\, T_2T_1 \times T_2,\, T_3 T_2 T_1 \times T_3 T_2 \times T_3,\dots,T_k \cdots T_1 \times \cdots \times T_k T_{k-1} \times T_k
\end{equation*}
is ergodic in the corresponding product space $(X^i,\mathscr{B}^i,\mu^i)$ then
\begin{equation*}
\lim_{N \to \infty} \frac{1}{\haar(\Phi_N)} \int\limits_{\Phi_N} \prod_{i=1}^k T_k^g \cdots T_i^g f_i \intd\haar(g) = \prod_{i=1}^k \int f_i \intd\mu
\end{equation*}
for any $f_1,\dots,f_k \in \lp^\infty(X,\mathscr{B},\mu)$.
\begin{proof}
It suffices to prove that the Furstenberg joining associated to the actions $T_1,\dots,T_k$ is the product measure $\mu^{k+1}$.
First note that when $k = 1$ this follows from Proposition~\ref{meanErgodic} because $\mathscr{C}_1$ is trivial by hypothesis.

Let $T_1,\dots,T_k$ be commuting, measurable actions satisfying the above ergodicity assumptions.
Assume by induction that $\nu_{k-1}$, the Furstenberg joining for the actions $T_1,\dots,T_{k-1}$, is the product measure $\mu^k$.
We know from Theorem \ref{UsingVanDerCorput} that
\begin{equation}
\label{afterVdcForJointErgodicity}
\int f_1 \otimes \cdots \otimes f_{k+1} \intd\nu_k = 0
\end{equation}
whenever $\condex{f_1 \otimes \cdots \otimes f_k}{\mathscr{I}_k} = 0$.
By hypothesis $T_k \cdots T_1 \times \cdots \times T_k$ is ergodic.
Therefore
\begin{equation*}
\condex{f_1 \otimes \cdots \otimes f_k}{\mathscr{I}_k} = \int f_1 \otimes \cdots \otimes f_k \intd\nu_{k-1} = \int f_1 \intd\mu \,\cdots \int f_k \intd\mu
\end{equation*}
for any $f_1,\dots,f_k$ in $\lp^\infty(X,\mathscr{B},\mu)$.
Thus \eqref{afterVdcForJointErgodicity} holds whenever $\int f_i \intd\mu = 0$ for some $1 \le i \le k$ as desired.
\end{proof}
\end{theorem}

One reason for being interested in Bergelson's conjecture is that it guarantees the existence of certain structures in large subsets of locally-compact, second-countable, amenable groups.
To make this precise, one needs to settle on a notion of largeness and then describe a correspondence principle that produces relevant measure-preserving actions from such sets.
This has been done in \cite{MR2561208}.
The correspondence principle does not yield \emph{ergodic} measure-preserving actions, so we cannot deduce combinatorial results from Theorem~\ref{thm:amenableFourTermSz}.
However, no ergodicity assumptions were made in the proof of Theorem~\ref{thm:amenableRoth}, and we now turn to combinatorial consequences of this result, discrete versions of which appear in \cite{MR1603193} and \cite{MR1481813}.

Given an invariant mean $M$ on $G$, we will say that a subset $S$ of $G$ is \define{substantial} if one can find a measurable subset $W$ of $G$ with $M(W) > 0$ and a symmetric open neighbourhood $U$ of 1 in $G$ such that $S \supset UW$.
Throughout this section we assume $G$ is infinite.

\begin{theorem}
\label{thm:correspondencePrinciple}
Given an invariant mean $M$ on $G$ and substantial subsets $S_1,\dots,S_k$ of $G$ one can find $c > 0$, a measurable action $T$ of $G$ on a separated, countably generated probability space $(X,\mathscr{B},\mu)$, and sets $B_1,\dots,B_k$ in $\mathscr{B}$ with positive measure such that
\begin{equation*}
M(g_1^{-1}S_1 \cap \cdots \cap g_k^{-1} S_k) \ge c\mu((T^{g_1})^{-1} B_1 \cap \cdots \cap (T^{g_k})^{-1} B_k)
\end{equation*}
for any $g_1,\dots,g_k$ in $G$.
\begin{proof}
The only discrepancies with Theorem 1.1 in \cite{MR2561208} are that $(X,\mathscr{B},\mu)$ is separated and countably generated, and that the action is measurable.
To overcome the first, note that since $G$ is second-countable the space $X$ obtained via the Gelfand representation in the proof of Theorem 1.1 in \cite{MR2561208} is a compact metric space.
Since the action obtained in \cite{MR2561208} is weakly measurable, using \cite{MR804490} we can assume the action is measurable.
\end{proof}
\end{theorem}

\begin{theorem}
Let $M$ be an invariant mean on $G \times G$ and let $S$ be a substantial subset of $G \times G$.
Then
\begin{equation*}
\{ g\in G \,:\, M(\{ (a,b) \in G \times G \,:\, (a,b), (a,gb), (ga, gb) \in S \}) > 0 \}
\end{equation*}
is syndetic.
\begin{proof}
The product $G \times G$ is also a locally-compact, second-countable, amenable group.
Let $S$ be a substantial subset of $G \times G$.
By the above theorem we can find an action $T$ of $G \times G$ on a separated, countably generated probability space $(X,\mathscr{B},\mu)$, a set $B \in \mathscr{B}$ having positive measure and some $c > 0$ such that
\begin{equation*}
M(S \cap (1,g)^{-1}S \cap (g,g)^{-1}S) \ge c \mu(B \cap (T^{(1,g)})^{-1} B \cap (T^{(g,g)})^{-1} B)
\end{equation*}
for every $g \in G$.
Define commuting actions $T_1$ and $T_2$ of $G$ on $(X,\mathscr{B},\mu)$ by $T_1^g = T^{(g,1)}$ and $T_2^g = T^{(1,g)}$.
The above becomes
\begin{equation}
\label{trianglesCorrespondence}
M(S \cap (1,g)^{-1}S \cap (g,g)^{-1}S) \ge c \mu(B \cap (T_2^g)^{-1}B \cap (T_2^g T_1^g)^{-1}B)
\end{equation}
for every $g \in G$.
By the discussion at the end of Section~\ref{sec:recurrenceResults}, the right-hand side of \eqref{trianglesCorrespondence} is positive for a syndetic set of $g \in G$.
\end{proof}
\end{theorem}

For our second result we need some facts about sets of recurrence.
A subset $S$ of $G$ is said to be \define{good for double recurrence} if the intersection of $S$ with
\begin{equation}
\label{eqn:rothRecurrenceTimes}
\{ g \in G \,:\, \mu(B \cap (T_2^g T_1^g)^{-1}B \cap (T_2^g)^{-1}B) > 0 \}
\end{equation}
contains an element different from the identity for any commuting, measurable actions $T_1,T_2$ of $G$ on a separated, countably generated probability space $(X,\mathscr{B},\mu)$ and any $B \in \mathscr{B}$ with positive measure.

\begin{lemma}
\label{lem:partitionRecurrence}
If a subset $S$ of $G$ good for double recurrence is finitely partitioned then one of the cells of the partition is good for double recurrence.
\begin{proof}
Fix a partition $S_1 \cup \cdots \cup S_r$ of a set $S$ good for double recurrence.
Suppose none of the $S_i$ is good for double recurrence.
Then for each $i$ one can find a separated, countably generated probability space $(X_i,\mathscr{B}_i,\mu_i)$ equipped with commuting measurable actions $T_{i,1}$ and $T_{i,2}$ such that
\begin{equation}
\label{eqn:partitionDoubleRecurrence}
\mu_i((T_{i,2}^gT_{i,1}^g)^{-1}B_i \cap (T_{i,2}^g)^{-1} B_i \cap B_i) = 0
\end{equation}
for all $g \ne 1$ in $S_i$.
Let $(X,\mathscr{B},\mu)$ be the product of the above probability spaces and let $B = B_1 \times \cdots \times B_r$.
Let $T_1$ and $T_2$ be the products of the $T_{1,i}$ and the $T_{2,i}$ respectively.
Since $S$ is good for double recurrence some $S_i$ contains an element $g$ different from the identity such that
\begin{equation*}
0 < \mu((T_2^g T_1^g)^{-1} B \cap (T_2^g)^{-1} B \cap B) = \prod_{i=1}^r \mu_i((T_{i,2}^g T_{i,1}^g)^{-1} B_i \cap (T_{i,2}^g)^{-1} B_i \cap B_i)
\end{equation*}
contradicting \eqref{eqn:partitionDoubleRecurrence}.
\end{proof}
\end{lemma}

\begin{lemma}
\label{lem:removingZeroDensityRecurrence}
If $S$ is a measurable subset of $G$ with $\upperdens(S) = 0$ then $G \backslash S$ is good for double recurrence.
\begin{proof}
First note that $\upperdens(G) \le \upperdens(G \backslash S) + \upperdens(S)$ so $\upperdens(G \backslash S) = 1$.
Passing to a sub F{\o}lner sequence we can assume that $\dens(G \backslash S) = 1$.
Since \eqref{eqn:rothRecurrenceTimes} has positive lower density with respect to this F{\o}lner sequence it cannot be disjoint from $G \backslash S$.
\end{proof}
\end{lemma}

Our second result concerns a non-commutative version of Schur's theorem that generalizes the discrete version in \cite{MR1603193}.
Let $Z(g)$ denote the centralizer of $g$ in $G$.
From \cite{MR1603193} we know that $\{ g \in G \,:\, [G:Z(g)] < \infty \}$ is a subgroup of $G$.
Moreover, it is measurable because it consists of those points in $G$ that have finite orbit under the action of $G$ on itself by conjugation and hence is a countable union of closed sets.
Given a subset $C$ of $G$ denote by $\tilde{C}$ the subset of $G \times G$ consisting of those $(a,b)$ such that $ab^{-1}$ belongs to $C$.

\begin{lemma}
\label{lem:densityOfTwiddleSet}
For any F\o{}lner sequence $\Psi$ in $G$ there is an increasing sequence $k_N$ such that
\begin{equation*}
\frac{\haar(\Psi_{k_N} \symdiff g \Psi_{k_N})}{\haar(\Psi_{k_N})} \le \frac{1}{N}
\end{equation*}
for all $g \in \Psi_N$. Moreover, if a measurable subset $S$ of $G$ has density with respect to $\Psi_N$ then $\tilde{S}$ has the same density with respect to $\Psi_N \times \Psi_{k_N}$.
\begin{proof}
The first part follows from the fact that $\haar(\Psi_N \symdiff g\Psi_N)/\haar(\Psi_N)$ converges to 0 uniformly on compact subsets of $G$. For the second part, if $\dens(S)$ exists then
\begin{align*}
\dens(S) = & \, \lim_{N \to \infty} \frac{\haar(S \cap \Psi_{k_N})}{\haar(\Psi_{k_N})}\\
 = & \, \lim_{N \to \infty} \frac{1}{\haar(\Psi_N)} \int 1_{\Psi_N}(g) \, \frac{\haar(S \cap g\Psi_{k_N})}{\haar(\Psi_{k_N})} \intd\haar(g)\\
 = & \, \lim_{N \to \infty} \frac{1}{\haar(\Psi_N) \haar(\Psi_{k_N})} \iint 1_{\tilde{S}}(g,h) 1_{\Psi_N}(g) 1_{\Psi_{k_N}}(h) \intd\haar(h) \intd\haar(g)
\end{align*}
as desired.
\end{proof}
\end{lemma}

\begin{theorem}
Suppose the subgroup $A = \{ g \in G \,:\, [G:Z(g)] < \infty \}$ of $G$ does not have finite index.
Then for any partition $C_1 \cup \cdots \cup C_r$ of $G$ and any open neighborhood $U$ of 1 in $G$ one can find $1 \le i \le r$ such that $UC_iU$ contains a subset of the form $\{ x,y,xy,yx \}$ with $xy \ne yx$.
\begin{proof}
Fix a partition $C_1 \cup \cdots \cup C_r$ of $G$ and an open neighborhood $U$ of 1 in $G$.
Let $V$ be an open neighborhood of 1 such that $VV \subset U$.
Permute the indices so that $VC_1V,\dots,VC_sV$ have positive upper density and $VC_{s+1}V,\dots,VC_rV$ have zero upper density with respect to $\Phi$.
Since $A$ has infinite index it has zero density.
Thus using Lemmas \ref{lem:partitionRecurrence} and \ref{lem:removingZeroDensityRecurrence} we can find $1 \le i \le s$ such that $VC_iV \backslash A$ is good for double recurrence.
Write $C = VC_iV$.
By passing to a sub F{\o}lner sequence we can assume that $\dens(C)$ exists and is positive.
By Lemma \ref{lem:densityOfTwiddleSet} we can find a F{\o}lner sequence $\Psi$ in $G \times G$ with respect to which the density of $\tilde{C}$ exists and is positive.
If $g \notin A$ then $Z(g)$ has infinite index and thus zero density, so $\tilde{Z}(g)$ will have zero density with respect to $\Psi$.

Let $M$ be a mean on $G \times G$ that agrees with $\dens_\Psi$ on the sets having density along $\Psi$.
The set $S = \{ (a,b) \,:\, ab^{-1} \in U C_i U \}$ contains the substantial subset $(V \times V)\tilde{C}$ of $G$.
Thus by Theorem~\ref{thm:correspondencePrinciple} we can find a measurable action $T$ of $G \times G$ on a separated, countably generated probability space $(X,\mathscr{B},\mu)$, some $B$ in $\mathscr{B}$ with $\mu(B) > 0$ and some $c > 0$ such that
\begin{equation}
\label{schurCorrespondence}
M(S \cap (1,g^{-1})S \cap (g^{-1},g^{-1})S) \ge c \int 1_B \cdot T^{(1,g)} 1_B \cdot T^{(g,g)} 1_B \intd\mu
\end{equation}
for any $g \in G$.
Putting $T_1^g = T^{(g,1)}$ and $T_2^g = T^{(1,g)}$ we can choose $g$ in $C \backslash A$ such that the right-hand side of \eqref{schurCorrespondence} is positive.
Since $M(\tilde{Z}(g)) = 0$ we can find $(a,b)$ in
\begin{equation*}
S \cap (1,g^{-1})S \cap (g^{-1},g^{-1})S \,\big\backslash\, \tilde{Z}(g)
\end{equation*}
giving $\{ ab^{-1}, ab^{-1}g^{-1}, g a b^{-1} g^{-1}, g \} \subset U C_i U$.
Putting $x = a b^{-1} g^{-1}$ and $y = g$ gives the desired result because $ab^{-1}$ does not belong to $Z(g)$.
\end{proof}
\end{theorem}

\printbibliography

\end{document}